\documentclass[11pt]{amsart}
\usepackage{amscd,amsmath,amssymb,amsthm}
\usepackage{color}
\usepackage{textcomp}
\usepackage{graphicx}
\usepackage{tikz}

\newtheorem{theorem}{{\bfseries Theorem}}[section]

\newtheorem{proposition}[theorem]{{ Proposition}}
\newtheorem{lemma}[theorem]{{ Lemma}}
\newtheorem{corollary}[theorem]{{ Corollary}}
\newtheorem{definition}[theorem]{{ Definition}}
\newtheorem{example}[theorem]{{ Example}}
\newtheorem{remark}[theorem]{{ Remark}}

\newcommand{\bt}{\begin{theorem}}
\newcommand{\et}{\end{theorem}}
\newcommand{\bl}{\begin{lemma}}
\newcommand{\el}{\end{lemma}}
\newcommand{\bp}{\begin{proposition}}
\newcommand{\ep}{\end{proposition}}
\newcommand{\bex}{\begin{example}}
\newcommand{\eex}{\end{example}}
\newcommand{\bc}{\begin{corollary}}
\newcommand{\ec}{\end{corollary}}
\newcommand{\bo}{\begin{proof}}
\newcommand{\eo}{\end{proof}}
\newcommand{\bd}{\begin{definition}}
\newcommand{\ed}{\end{definition}}
\newcommand{\br}{\begin{remark}}
\newcommand{\er}{\end{remark}}
\newcommand{\be}{\begin{enumerate}}
\newcommand{\ee}{\end{enumerate}}

\newcommand{\cM}{{\mathcal{ M}}}
 
\newcommand{\cO}{{\mathcal{ O}}}

\newcommand{\Z}{{\mathbb Z}}

\newcommand{\cS}{{\mathcal{ S}}}
\newcommand{\N}{{\mathbb N}}
\newcommand{\R}{{\mathbb R}}
\newcommand{\C}{{\mathbb C}}

\begin{document}

\title{Some relations in topological dynamics}

\author{Joseph Auslander}
\address{Department of Mathematics, University of Maryland,  College Park, MD 20742, USA.}

\email{ jna@umd.edu}

\author{Anima Nagar}
\address{Department of Mathematics, Indian Institute of Technology Delhi,
	Hauz Khas, New Delhi 110016, INDIA.}

\email{ anima@maths.iitd.ac.in}


\date{24 June, 2024}

\vspace{0.5cm}



	\begin{abstract}
Relations always play an important role in the study of topological dynamics. Proximal, distal and almost periodic relations are well studied in literature. We further  this direction and analogously study the strongly proximal and weakly distal relations. 

This gives a new class of flows - the weakly distal flows. We observe that the well known Morse-Thue substitution flows and Chacon transformations are weakly distal.
	
\end{abstract}

\maketitle

\section{Introduction}

A \emph{flow} $ (X, T) $ is a jointly continuous (left) action of the \emph{topological group}
$ T $ on the \emph{compact metric space} $ (X,d) $ or a compact Hausdorff space $(X, \mathfrak{T})$. A closed subset $A\subseteq X$ is \emph{invariant} if $TA :=\{ta\ :\ a\in A, \ t\in T\}\subseteq A$. If $A$ is \emph{invariant} then $ T\times A \rightarrow A$ is an action and $(A, T)$ is a flow,  called a \emph{subsystem} of $(X,T)$. Note that for all $x \in X$, the \emph{orbit closure} $ \overline{Tx} $ is an invariant set.

A point $x \in X$ is called \emph{almost periodic} if its orbit closure $ \overline{Tx} $ is a minimal set, i.e. is non-empty, closed, $ T- $invariant, and minimal with respect to these
properties. If $(X,T)$ is a \emph{minimal flow} then all $x \in X$ are almost periodic.

\bigskip 

We refer to \cite{AUS} for the general theory.

\bigskip

Each $t \in T$ is identified with the map $x \to tx, \ x \in X$. Thus $T \subset X^X$, endowed with the product topology, and its \emph{enveloping semigroup}  is $E = E(X) = E(X,T) := \overline{T}$ in $X^X$. We recall the  connections between dynamical properties of the flow $ (X,T) $ and algebraic properties of its enveloping semigroup $ E(X)$. Also $E(X)$ contains idempotents, i.e. there exists $u \in E(X)$ with $u^2 = u$.

Let $ M $ denote an arbitrary  \emph{minimal left-ideal} in $ E(X) $, i.e. $E(X)M \subseteq M$. Henceforth,  a left-ideal will be just called \emph{ideal}. $ J_M $ denotes the collection of \emph{idempotents} in the minimal ideal $M$ - these are called \emph{minimal idempotents} and $J$ denotes the collection of all minimal idempotents in $E(X)$. When $u \in J_M$ then $uM$ is a group with identity $u$. Also $Mp = M$ for all $p \in M$ and $pu = p$ for all $p \in M$, \ $u \in J_M$. And for $u \neq v \in J_M$, $uM \cap vM = \emptyset$. 

 Consider the product flow $(X^n,T)$, with the diagonal action  defined as $(x_1, x_2, \ldots, x_n) \stackrel{t}{\to} (tx_1, tx_2, \ldots, tx_n)$. Since the enveloping semigroup is endowed with the product topology i.e. the topology of pointwise convergence of functions, the enveloping semigroup $E(X^n) = \Delta(E(X)^n)$ - the diagonal in the product $E(X)^n$. For  $p \in E(X)$, we shall denote the corresponding tuple in $E(X^n)$ also by $p$.

 It is well known that $\overline{T(x)}$ is minimal $\Longleftrightarrow$ for every minimal ideal $ M  $ in $E(X)$, we have $\overline{T(x)} =  Mx$ $\Longleftrightarrow$  there is a $v \in J_M$ such that $ vx = x $.

\bigskip

For flows $(X,T)$ and $(Y,T)$, the continuous surjection $\pi: X \to Y$ defines a \emph{factor map ( also called homomorphism or extension)} if $\pi(tx) = t\pi(x)$ for all $x \in X$ and $t \in T$, i.e. $\pi$ is continuous, surjective and equivariant. $Y$ is then called  a \emph{factor} of $X$ and $X$ is called an \emph{extension} of $Y$.  Further when $\pi$ is a homeomorphism then it is called a \emph{topological conjugacy} and the two flows $(X,T)$ and $(Y,T)$ are \emph{conjugate}, i.e. display the same dynamical properties.

The factor map $\pi: X \to Y$ induces the factor map $\theta: E(X) \to E(Y)$. Note that $ \pi(px)=\theta(p) \pi(x) $. Let $ \mathfrak{M}  $ be the collection of all distinct minimal ideals in $E(X)$, then $\{\theta (M) : M \in \mathfrak{M}\} $ is the collection of all distinct minimal ideals in $E(Y)$. 

Let  $M,M'$ be distinct minimal ideals in $E(X)$. We call  $u \in J_M$ and $u' \in J_{M'}$ to be equivalent  if $uu' = u'$ and $u'u = u$ and write $u \sim u'$. Now for each $u \in J_M$, there exists $u' \in J_{M'}$ such that $u \sim u'$. If $\{u_M : M \in \mathfrak{M}\} $ is a  collection of  equivalent minimal idempotents in $E(X)$, then $\{\theta (u_M) : M \in \mathfrak{M}\} $ is the collection of  equivalent minimal idempotents in  $E(Y)$.

Also note that for  $u,v \in J_M$, we have $u \nsim v$. Since $ uv = u$ being in $M$ and if $u \sim v$ then $uv = v$ giving $u = v$.

 From \cite{AUS}  recall, for a  flow $(X,T)$ a set $A\subset X$ is said to be an \emph{almost periodic set} if any point $z\in X^{|A|}$ with
$range(z)=A$ is an almost periodic point of $(X^{|A|},T)$, where $|A|$ is the cardinality of $A$. 

 For each $u \in J$, consider the set $F_{u}=\lbrace x\in X: ux=x \rbrace = uX$. We note that if $u \sim u'  (\in J)$, then $F_u = F_{u'}$ and these are almost periodic sets(c.f. \cite{An}). 

\begin{lemma} \cite{An} \label{aps} For the flow $(X,T)$, the following are equivalent:
	
	1. $A \subset X$ is an almost periodic set.
	
	2. There exists a minimal idempotent $u \in E(X)$ such that $ua = a$,  for every $ a \in A$. \end{lemma}

The flow $(X,T)$ induces the flow $(2^X,T)$, where $2^X$ is the space of all closed subsets of $X$ endowed with the Hausdorff topology and the action of $T$ given as $t(A) =\{ta: a \in A\}$, for all $A \in 2^X$.

\section{Relations and Elementary Properties}

Let $\Delta$ denote the \emph{diagonal} in $X \times X$. 

\subsection{Almost Periodic, Proximal and Distal   Relations} We recall some known relations and study their properties.

\bigskip

$ \bigstar $ Let $\Omega (= \Omega(X))$ denote the set  of all \emph{almost periodic points} of
the product flow $ (X \times X, T) $. Then, $\Omega$ can be considered as a relation on $X$. Thus, the \emph{almost periodic relation} $\Omega$ is defined as
$$\Omega(X) = \Omega := \{(x,y) \in X \times X: \ \text{there exists a} \  v \in J \ \text{such that} \ v(x,y) = (vx,vy) = (x,y) \}.$$ 

$\Omega$ is a symmetric, $ T- $ invariant
relation. If $(X,T)$ is minimal then $\Omega$ is also reflexive though   is not always an equivalence relation.  Also if $(x,y) \in \Omega$ then both $x,y \in X$ are almost periodic points. 

For almost periodic $x \in X$, define
$$\Omega[x] := \{y \in X : (x,y) \in \Omega\}$$ 

-  the set of all  almost periodic points whose product with $x \in X$ is also almost periodic.

\begin{lemma}
	For a flow $(X,T)$, $\Omega[x] = \bigcup \limits_{u \in J: ux = x} uX $, for all almost periodic $x \in X$.
\end{lemma}

\begin{proof} 
	Let $x \in X$ be almost periodic. For $y \in X$ if $(x,y) \in \Omega$ then there exists $u \in J$ such that $u(x,y) = (x,y)$. Thus, $ux = x$ and $y \in uX$.
	
	Conversely, let $u \in J$ such that $ux = x $ and $ \ y \in uX$. Then $(x,y) \in \Omega$.
\end{proof}

\bigskip 

$ \bigstar $ The \emph{proximal  relation} is defined as

$$P(X) = P := \{(x,y) \in X \times X: \overline{T(x, y)} \cap \Delta \neq \emptyset\}. $$

$P$ is reflexive, symmetric, $ T-$invariant
relation, and in general is not an equivalence
relation.

The  flow $ (X,T) $ is called \emph{proximal} if and only if the relation $P =  X \times X$.

The set of all points proximal  to $x \in X$ is
$$P[x] := \{y \in X : (x,y) \in P\}.$$ 

It is known that if $ (X,T) $ is minimal, then  $	P[x] = \{vx: v \in J\}$.

Note that the pair $ (x,y) \in P $ $\Longleftrightarrow$ there exists a minimal ideal $ M $ in $ E(X) $ such that $ px = py $ for all $p \in M$. And  $(x,vx) \in P(X)$ for all $v \in J$, $x \in X$. 
For every $x \in X$, there exists an almost periodic point $x' \in X$ such that $(x,x') \in P$.

We use this trivial yet important lemma.

\begin{lemma}
	For a flow $(X,T)$ there are no non-trivial proximal pairs which are almost periodic in the product flow $ (X \times X, T )$, i.e. $P \cap \Omega \subset \Delta$.
\end{lemma}

\begin{proof} Let $(x,y) \in P$ then there exists a net $\{t_\alpha\}$ in $T$ with $t_\alpha(x,y) \to (z,z)$ for some $z \in X$. Now since $(x,y) \in \Omega$ there exists a net $\{s_\beta\}$ in $T$ with $s_\beta (z,z) \to (x,y)$ and so $x=y$.
	
	Thus $(x,y) \in P \cap \Omega \ \Longrightarrow (x,y) \in \Delta$.
\end{proof}

Note some important characterizations here:

\begin{theorem} \cite{AUS} \label{unique}
	Let $ (X, T) $ be a 
flow. Then the following are equivalent:

1. $  P $ is an equivalence relation 

2. $ E(X) $ contains a unique minimal ideal.

3.  $ (x, y) \in P \Longrightarrow  \overline{T(x, y)} \subset P $.
\end{theorem}

\begin{lemma}\cite{AUS} 
	Let $ (X, T) $ be a 	flow. Then $(x,ux) \in P$ for $u \in J$.
\end{lemma}

\bigskip

$ \bigstar $ A point $x \in X$ is called \emph{distal} if for every $y \neq x \in X$ we have $(x,y) \notin P$, i.e. $P[x] = \{x\}$. It is known that a distal point is almost periodic.

\bigskip

A pair  $x,y \in X$ is called \emph{distal} if $(x,y) \notin P$. The \emph{distal relation} is defined as

$$D(X) = D := \{(x,y) \in X \times X: \overline{T(x, y)} \cap \Delta = \emptyset\}. $$

$D$ is  symmetric, $ T-$invariant
relation, and cannot be reflexive since $D \cap \Delta = \emptyset$.

Thus, $(x,y) \in P \ \Longleftrightarrow \ (x,y) \notin D \ \text{and} \ (x,y) \in D \ \Longleftrightarrow \ (x,y) \notin P$. And $X \times X = P \cup D$ is a disjoint union.

The set of all points distal  to $x \in X$ is

$$D[x] = \{y \in X : (x,y) \in D\}.$$

For all $x \in X$, $X = P[x] \cup D[x]$ is a disjoint union. Note that $D[x] \supset \bigcap \limits_{u \in J} uX \setminus \{x\} $ - the set of all distal points in $X \setminus \{x\}$. 
Note that $ x \in X $ is a distal point $\Longleftrightarrow$  $ ux = x $ for all $  u \in J $. Every pair of distinct points in $ vX = \{x: vx = x\} $ is distal for all $v \in J$, i.e. $ \bigcup \limits_{x \neq y \in vX} (x,y) \subset D$ for all $v \in J$.

\bigskip

We say that the relation $D$ is \emph{full} if $D = X \times X \setminus \Delta$, i.e., $P= \Delta$.

\bigskip

	 The  flow $ (X,T) $ is called \emph{{distal}} if and only if the relation $D$ is full, i.e., $D = X \times X \setminus \Delta$, i.e., $P= \Delta$.  
	 Distality can be characterized as:

\begin{theorem} \cite{AUS} The flow  $ (X,T) $ is distal $\Longleftrightarrow$ $ E(X) $ is a group $\Longleftrightarrow$  $X= \bigcap \limits_{v \in J} vX$ $\Longleftrightarrow$ $P=\Delta$. \end{theorem}

Note that $\overline{T(x,y)} \cap P \neq \emptyset \Longrightarrow (x,y) \in P$.

\begin{corollary} 
	For a flow $ (X, T) $,  $(x,y) \in D \ \Longleftrightarrow \ \overline{T(x, y)} \subset D $.
\end{corollary}

\begin{proof} 
	We note that $D$ is $T-$invariant. Let $(x,y) \in D $. Suppose $(a,b) \in  \overline{T(x, y)}$, then there exists $p \in E(X \times X)$ such that  $p(x,y) = (a,b)$. Thus, $ \overline{T(a, b)} \subset \overline{T(x, y)}$ and so $ \overline{T(a, b)} \cap  \Delta = \emptyset $. Thus $\overline{T(x, y)} \subset D$.
	
	The converse is obvious.
\end{proof}

\bigskip

 \begin{lemma}
 	For a flow $(X,T)$,	$  \Omega = X \times X \Longleftrightarrow (X,T) $ is distal $  \Longleftrightarrow D = X \times X \setminus \Delta \Longleftrightarrow P = \Delta$.\end{lemma}
 \begin{proof}
 	We note that if $(x,y) \in \Omega$ then there exists $u \in J$ such that $u(x,y) = (x,y)$, i.e. $x,y \in uX$. Now for given $x \in X$, we have $D[x] \supset uX \setminus \{x\}$ $\Longleftrightarrow (x,y) \in \Omega$( with $x \neq y$)  $\Longleftrightarrow \ (x,y) \in D$ $\Longleftrightarrow D$ is full $ \Longleftrightarrow P = \Delta$.
 \end{proof}

\begin{remark}
	Let $(Z,T)$ be a subsystem of $(X,T)$, then $P(Z) \subset P(X)$, $D(Z) \subset D(X)$ and $\Omega(Z) \subset \Omega(X)$.
	
\end{remark}

\bigskip

For flows $(X,T)$ and $(Y,T)$ let $\pi: X \to Y$ define a factor map. Then
$$R_\pi = R_\pi(X) = \{(x,x') \in X \times X: \pi(x) = \pi(x')\}$$
 is a \emph{closed $T-$invariant equivalence relation (icer)}.
 
 Every icer $ R $ on
 $ X $ determines a factor $ X/R $. If $ M $ is a minimal set in $ X $, then $ \pi(M) $ is minimal in $ Y $.

\begin{definition}
	For flows $(X,T)$ and $(Y,T)$ let $\pi: X \to Y$ define a factor map. The  map $\pi$ is called \emph{proximal ( resp. distal)} when for every $y \in Y$ if $x_1,x_2 \in \pi^{-1}(y)$ then $(x_1,x_2) \in P(X)(\text{resp.} \ D(X))$. 

 This factor map $\pi$ is called \emph{highly proximal (HP)} if every fiber can be shrunk uniformly to a point, that is given $ y \in  Y $ there is a net $ \{t_n\} $ in $ T $ such that   $ t_n \pi^{-1}(y) \to \{x\} $, for some $ x  \in X$,  wrt the Hausdorff topology in $2^X$.
\end{definition} 

A highly proximal factor is proximal but not conversely. 
It is known that if $ X $ and $ Y $ are metric spaces, then $\pi$ is highly proximal if and only if it is almost one-one, i.e.,  there is a residual subset of $ Y $ on which $ \pi $ maps one to one. A  factor map is an isomorphism if and only if it is open and highly proximal.

\bigskip

 \begin{lemma} \label{rwy=ry}
 	For a flow $ (Y,T) $,	let  $w \in E(Y) $  be such that $ (wy,y) \in P(Y) $ for all $ y \in Y $.  Then $ w^2=w $.
 \end{lemma}
 
 \begin{proof}
 	Observe that our hypothesis implies that $(w(wy),wy) \in P(Y)$ for all $ y \in Y $. But $ (w(wy),wy) \in \Omega(Y) $ and so $ wwy = wy $. Since $ y $ is arbitrary, we have
 	$ w^2=w $.
 	
 \end{proof}
 
 Suppose $w \in E(Y)$ is a minimal idempotent, and let $w = \psi(u)$ for minimal idempotent $u \in \beta T$.  A natural question arises: \emph{What are all the elements in $\theta^{-1}(w)$?  }

 \begin{theorem}
 	Let the factor $ \pi:X \to Y $, with $ (Y,T) $ minimal, induce the factor  $ \theta:E(X) \to E(Y) $. Let $ K $ be a minimal  ideal in $ E(X) $ and
 	let $ p \in K  $ with $ w=\theta(p) $. Then $ w^2=w $ if and only if   $ \pi(px,x ) \in P(Y)$ for every $ x \in X $.
 	
 \end{theorem}
 
 \begin{proof}
 	Suppose $ w^2=w $. Then  $\pi(px,x) = (w \pi(x),\pi(x)) \in P(Y)$ by Lemma \ref{rwy=ry}.
 	
 	Conversely  we have $ (\theta(p)\pi(x),\pi(x)) \in P(Y)$ for every $ x \in X $. By Lemma \ref{rwy=ry}, it follows that $w^2 = w$.
 \end{proof}

 Thus for a minimal idempotent $w \in E(Y)$, we have 
 $$\theta^{-1}(w) = \bigcup \limits_{\{K \ \text{minimal ideal in} \ E(X): \ w \ \in \ \theta(K)\}} \ \{p \in K : \ \text{for every} \  x \in X, \  \pi(px,x) \in P(Y) \}.$$ 
 
 \bigskip

\begin{theorem} \label{pi by pi}
	For flows $(X,T)$ and $(Y,T)$ and the factor map $\pi: X \to Y$ we have:

\begin{enumerate} 
	 
	\item $\pi \times \pi(P(X)) \subset P(Y)$ 
	\item  $\pi \times \pi(D(X)) \supset D(Y)$ 
	\item $\pi \times \pi(\Omega(X)) = \Omega(Y)$.
\end{enumerate}
\end{theorem} 

\begin{proof}
	Note that $ \pi \times \pi (P(X)) \subset P(Y) $ is clear. And by complementation we have $ \pi \times \pi (D(X)) \supset D(Y) $. It is clear that $ \pi \times \pi (\Omega(X)) \subset \Omega(Y) $. 	Let $(y,y') \in \Omega(Y)$, then there is a minimal idempotent $u \in E(X)$ such that for the induced homomorphism $\theta u (y,y') = (y,y')$. Let $ \pi \times \pi(x, x') = (y, y') $. 	Then $ \pi \times \pi(ux, ux') = (\theta u (y), \theta u (y')) = (y, y') $ and  so $ (ux, ux') \in \Omega(X) $.
\end{proof}

\begin{theorem} \label{pi by pi inverse}
	For flows $(X,T)$ and $(Y,T)$ let $\pi: X \to Y$ define a factor map. Then
	\begin{enumerate}
		
		\item $P(X) \subset (\pi \times \pi)^{-1} P(Y)$ with equality when $\pi$ is proximal. 
		\item  $D(X) \supset (\pi \times \pi)^{-1} D(Y)$ with equality when $\pi$ is distal.
		\item $\Omega(X) \subset (\pi \times \pi)^{-1} \Omega(Y)$ with equality when $\pi$ is distal. 
	\end{enumerate}
	
\end{theorem}
\begin{proof}
	Note that, we always have 
$ \pi \times \pi(P(X)) \subset P(Y)$, 
   $ \pi \times \pi(D(X)) \supset D(Y)$ and $ \pi \times \pi(\Omega(X)) \subset  \Omega(Y)$ by Theorem \ref{pi by pi}.
 
 \begin{enumerate}
 	\item Let $\pi$ be proximal. Now  since $ P(X) \subset (\pi \times \pi)^{-1}(P(Y)) $ we prove the converse.  Let $ (y, y') \in P(Y) $  
 	and $\pi \times \pi(x, x') = (y, y') $,  and let $ u \in E(X) $ be a minimal idempotent for which $ \theta u(y) = \theta u(y') $. Then $ (ux, ux') \in \Omega(X)  $ and $ \pi(ux) = \pi(ux') $. Since $ \pi $ is proximal, $ (ux, ux') \in P(X)$. And so $ ux =
 	ux' $ and $ (x, x') \in P(X) $.
 	
 	\item Let $\pi$ be distal. Now $ D(X) \supset (\pi \times \pi)^{-1}(D(Y)) $ and we prove the converse.  Let $ (x, x') \in D(X) $  
 	and $\pi \times \pi(x, x') = (y, y') $.   Suppose $ u \in E(X) $ is a minimal idempotent for which $ \theta u (y) = \theta u(y') $. Then  $ \pi(ux) = \pi(ux') $ and since $ \pi $ is distal, we have $ ux =
 	ux' $ contradicting that $ (x, x') \in D(X) $. So $(y,y') \in D(Y)$.
 	
 	\item Assume that $\pi$ is distal. Suppose that $ (y, y') \in \Omega(Y) $ and
 	$ \pi \times \pi(x, x') = (y, y') $. Let $ u  \in E(X) $ be a minimal idempotent  such that $ (\theta u(y), \theta u(y')) = (y, y') $. 
 	
 	Then $ \pi \times \pi(ux, ux') = (\theta u(y), \theta u(y')) = (y, y') $ and $ (ux, ux') \in \Omega(X) $. Further, $ (ux, x) \in P(X ) $ and $ \pi(ux) = \theta u (y) = y = \pi(x) $. Since $ \pi $  is distal, $ ux = x $. Similarly $ ux' = x' $. Therefore $ (x, x') = u(x, x') \in  \Omega(X) $.
 \end{enumerate}

\end{proof}

\bigskip

We can generalize these results for arbitrary products.

Let $(\{(X_\alpha, T): \alpha \in \Lambda\})$ be a family of flows with action by the same group $T$.

Let $X = \prod \limits_{\alpha \in \Lambda} X_\alpha$. Then $\pi_\alpha: X \to X_\alpha$ gives the projection on the $\alpha$th coordinate and $\theta_\alpha: E(X) \to E(X_\alpha)$ the induced homomorphism on enveloping semigroups, for each $\alpha$.

We recall this interesting result:

\begin{theorem} \cite{A}
	Let $(\{(X_\alpha, T): \alpha \in \Lambda\})$ be a family of flows with action by the same group $T$, and let $X = \prod \limits_{\alpha \in \Lambda} X_\alpha$. 
	
	Then	$ P(X) $ is an equivalence relation on $ X $ if and
	only if $ P(X_\alpha) $ is an equivalence relation on each $X_\alpha$.
\end{theorem}

Let $\{I_\nu\} $ be the collection of all distinct minimal ideals in $E(X)$, then $\{\theta_\alpha (I_\nu) = I_{\alpha,\nu}\} $ is the collection of all distinct minimal ideals in $E(X_\alpha)$. 

Now for each $u_\nu \in I_\nu$, there exists $u_{\nu'} \in I_{\nu'}$ such that $u_\nu \sim u_{\nu'}$. Let $\{u_\nu\} $ be a collection of  equivalent minimal idempotents in $I_\nu$ for each $\nu$, then $\{\theta_\alpha (u_\nu) = u_{\alpha,\nu}\} $ is a collection of  equivalent minimal idempotents in $I_{\alpha,\nu}$ for each $\nu$ and $\alpha$.

Let $(x,y) \in P(X)$, then there is a minimal ideal $I_\mu \subset E(X)$ such that $px =py$ for all $p \in I_\mu$. If $x = (x_\alpha), \ y = (y_\alpha)$, then $p_\alpha x_\alpha = p_\alpha y_\alpha$ for all $p_\alpha \in I_{\alpha,\mu}$ a minimal ideal in $E(X_\alpha)$.

Thus we have, let $(\{(X_\alpha, T): \alpha \in \Lambda\})$ be a family of flows with action by the same group $T$, and let $X = \prod \limits_{\alpha \in \Lambda} X_\alpha$. 	Then $\pi_\alpha \times \pi_\alpha P(X) \subset  P(X_\alpha)$ for every $\alpha \in \Lambda$.

\begin{remark}
	We note that equality in the expression in the Theorem above need not hold. Let $(X,T)$ be such that $P(X)$ is not an equivalence relation. Let $x,y,z \in X$ be such that $(x,y) , (y,z) \in P(X)$ but $(x,z) \notin P(X)$. If $((x,y), (y,z)) \in P(X \times X)$, then we have a $p \in E(X)$ such that $p(x,y) = p(y,z) = (a,b)$(say). But then $a = py$ and also $py = b$ gives $a = b$. Hence $px = a = b = pz$ gives $(x,z) \in P(X)$ - which is a contradiction.
	
	So $((x,y), (y,z)) \notin P(X \times X)$.
\end{remark}

Note that $\pi_\alpha( \overline{Tx}) \subset \overline{Tx_\alpha} = \overline{\pi_\alpha(Tx)}$ and so $\pi_\alpha (\overline{Tx}) = \overline{Tx_\alpha}$.

Thus $(x,y) \in D(X) $ implies that $\overline{T(x,y)} \cap \Delta = \emptyset$ which gives  $\overline{T(x_\alpha,y_\alpha)} \cap \Delta = \emptyset$ for some $\alpha$. Conversely, if  $\overline{T(x_\alpha,y_\alpha)} \cap \Delta = \emptyset$ for some $\alpha$ then $\overline{T(x,y)} \cap \Delta = \emptyset$. And so we have,

\begin{theorem} 
	Let $(\{(X_\alpha, T): \alpha \in \Lambda\})$ be a family of flows with action by the same group $T$, and let $X = \prod \limits_{\alpha \in \Lambda} X_\alpha$. 
	
	Then	$ (x,y) \in D(X) $  if and
	only if $(x_\alpha, y_\alpha) \in D(X_\alpha) $ for some $\alpha$.
\end{theorem}

Let $ (x,y) \in \Omega(X) $. Then there exists $u \in J \subset E(X)$ such that $u(x,y) = (x,y)$. But this implies $\theta_\alpha(u)(x_\alpha, y_\alpha) = (x_\alpha, y_\alpha)$, i.e., $(x_\alpha, y_\alpha)  \in \Omega(X_\alpha)$ for all $\alpha \in \Lambda$. Again if $(x_{\alpha_0}, y_{\alpha_0})  \in \Omega(X_{\alpha_0})$ for some $\alpha_0$ then as observed in Theorem \ref{pi by pi}, we get $(x_\alpha, y_\alpha)  \in \Omega(X_\alpha)$ for all other $\alpha \neq \alpha_0$ such that for $x = (x_\alpha), \ y = (y_\alpha)$ we have $\pi_{\alpha_0} \times \pi_{\alpha_0} (x,y) = (x_{\alpha_0}, y_{\alpha_0})$. Thus,

\begin{theorem}
	Let $(\{(X_\alpha, T): \alpha \in \Lambda\})$ be a family of flows with action by the same group $T$, and let $X = \prod \limits_{\alpha \in \Lambda} X_\alpha$. 
	
	Then $\pi_\alpha \times \pi_\alpha \Omega(X) =  \Omega(X_\alpha)$ for all $ \alpha \in \Lambda $.
\end{theorem}

\bigskip

We recall a characterization of distal extensions:

\begin{theorem} \cite{JoAu}
	Let $(X,T)$ and $(Y,T)$ be flows with the same acting group $T$, and let $ (Y,T) $ be minimal. For the extension $\pi: X \to Y$ the following are equivalent:
	
	(1) $\pi$ is distal.
	
	(2) Every fiber $\pi^{-1}(y)$ is an almost periodic set.
\end{theorem}

However, in general we have

\begin{theorem} \label{fiber}
	Let $(X,T)$ and $(Y,T)$ be flows with the same acting group $T$, and let $ (Y,T) $ be minimal. Then for the  extension $\pi: X \to Y$ every fiber $\pi^{-1}(y)$ contains an almost periodic set.
\end{theorem}

\begin{proof}
	Let $\pi: X \to Y$ be an extension, and let $y \in Y$. Since $y$ is almost periodic there exists a minimal idempotent $w \in E(Y)$ such that $wy=y$. Let $u \in E(X)$ be a minimal idempotent such that $\theta(u) = w$. Then the almost periodic set $u \pi^{-1}(y) \subset \pi^{-1}(y)$ (see Lemma \ref{aps}).
\end{proof}

\begin{remark}
	Note that if $u \nsim v$ are minimal idempotents in $E(X)$ then for any $A \subset X$  the almost periodic sets $uA \neq vA$, though it is possible that $uA \cap vA \neq \emptyset$. However, as discussed in \cite{An}, if $u \sim u'$ then the almost periodic sets $u'(uA) = u'A$.
	
	If  $\pi$ is distal then $\pi^{-1}(y)$ is itself an almost periodic set, and for some minimal idempotent $u \in E(X)$, following Lemma \ref{aps} we have $u \pi^{-1}(y) = \pi^{-1}(y)$.
	
	However if $\pi$ is not distal, there can be multiple minimal idempotents in $\theta^{-1}(w)$, and so $\pi^{-1}(y)$ may contain more than one almost periodic sets.
\end{remark}

	\bigskip
	
\begin{remark}	
	Consider $x, z \in \pi^{-1}(y)$, and $u,v \in J_I \cap \theta^{-1}(w)$ where $wy =y$ and $I \subset E(X)$ is a minimal ideal in $E(X)$. We note that if $u' \in J_{I'} \cap \theta^{-1}(w)$ such that $u' \sim u$, then $u'\pi^{-1}(y)$ need not be the same as $ u\pi^{-1}(y)$ and so while we look into the results in case of choosen $u,v \in J_I$ we will have  the same conclusions for $u', v' \in J_{I'}$.

	\begin{enumerate}
		\item Let $x =vz$ then $vx = vvz = vz$ and so $(x,z) \in P$.
		
		\item If $ux = z$ then $uux = ux =uz$ and so $(x,z) \in P$.
		
		\item For all $ x \in \pi^{-1}(y)$, $(x,ux) \in P$.
		
		\item For $x, z \in u\pi^{-1}(y)$ we have $(x,z) \in D$. Thus $[u \pi^{-1}(y) \times u \pi^{-1}(y) ] \setminus \Delta \subset D$.

		\item For $x,z \in \pi^{-1}(y)$, let $x \neq vz$. Then $u(ux,vz) = (ux, u(vz)) \in u\pi^{-1}(y) \subset D$ and so $u(ux) \neq u(vz) $.
		
		Now let $I^*$ be a minimal ideal in $E(X)$ such that $I^* \neq I$, and let $u^* \in J_{I^*}$ be such that $u^* \sim u$. Then $u^*(ux,vz) = (u^*(ux), u^*(uvz)) = (u^*x, u^*(vz)) \in  u^*\pi^{-1}(y) \subset D$ and so $u^*(ux) \neq u^*(vz)$.
		
		If $(ux,vz) \in P$ then there  exists a minimal ideal $K$ in $E(X)$ for which $p(ux) = p(vz)$ for all $p \in K$. But for $u_p \in K$ such that $u_p \sim u$, we have $u_p(ux) \neq u_p(vz)$ and so $(ux,vz) \in D$.
		
		Thus  $[u\pi^{-1}(y) \times v\pi^{-1}(y)] \setminus\Delta \subset D$.
		
		\item For $x,z \in \pi^{-1}(y)$, let $x \neq vz$. Then $v(x,vz) = (vx, vz) \in v\pi^{-1}(y) \subset D$ and so $v(x) \neq v(vz) $.

		Then as  observed above we get
		$[\pi^{-1}(y) \times v\pi^{-1}(y)] \setminus\Delta \subset D$.
		
		\item For $x,z \in \pi^{-1}(y)$, let $ux \neq z$. Then $u(ux,z) = (ux, uz) \in u\pi^{-1}(y) \subset D$ and so $u(ux) \neq u(z) $.
		
		Then as observed above we get 
		$[u\pi^{-1}(y) \times \pi^{-1}(y)] \setminus\Delta \subset D$.
		
		\item For $x,z \in \pi^{-1}(y)$, let $x \neq v'z, \ ux \neq z$ for $u \in I, \ v' \in I'$ with $I \neq I'$. Let $v \in I, u' \in I'$ be such that $u \sim u'$ and $v \sim v'$. Then $u'(ux,v'z) = (u'x, u'(v'z)) \in u'\pi^{-1}(y) \subset D$ and so $u'(ux) \neq u'(v'z) $. Similarly, $v(ux) \neq vz $.
		
		Now let $I^*$ be a minimal ideal in $E(X)$ such that $I^* \neq I, I'$, and let $u^* \in J_{I^*}$ be such that $u^* \sim u$. Then $u^*(ux,vz) = (u^*(ux), u^*(uvz)) = (u^*x, u^*(vz)) \in  u^*\pi^{-1}(y) \subset D$ and so $u^*(ux) \neq u^*(vz)$. Similarly, $v^*(ux) \neq v^*(z)$.
		
		If $(ux,vz) \in P$ then there  exists a minimal ideal $K$ in $E(X)$ for which $p(ux) = p(vz)$ for all $p \in K$. But for $u_p \in K$ such that $u_p \sim u$, we have $u_p(ux) \neq u_p(vz)$ and so $(ux,vz) \in D$.
		
		Thus  $[u\pi^{-1}(y) \times v'\pi^{-1}(y)] \setminus\Delta \subset D$.
		
		Similarly, $[u'\pi^{-1}(y) \times v\pi^{-1}(y)] \setminus\Delta \subset D$.
		
	\end{enumerate}
	
	\textbf{Question:} What can be said about the pair of points in $\pi^{-1}(y) \ \setminus \ \bigcup \limits_{u \in (J \cap \theta^{-1}(w))}	 u \pi^{-1}(y) $?
	
\end{remark}

\bigskip

\subsection{Strongly Proximal and Weakly Distal Relations}  We now define and discuss some new relations.

\begin{definition}
	
The \emph{strongly proximal  relation} is defined as 

$$SP(X) = SP= \{(x,y) \in X \times X: \overline{T(x, y)} \cap \Omega \subset \Delta \}. $$
\end{definition}

Clearly, $SP \subseteq P$.

Note that for any $x \in X$, and a minimal set $F \subset \overline{Tx}$,  then  $K = Ix$, for minimal ideals $I$ in $E(X)$. Hence we can equivalently say that,

\begin{theorem} For the flow $(X,T)$
	
$$(x,y) \in SP \ \Longleftrightarrow \  I(x,y) \subset \Delta, \ \text{for every minimal ideal $I$ in $E(X)$}.$$
\end{theorem}

\begin{remark}
	Let $ (X, T) $ be  a proximal flow, then $P = X \times X = SP$. When $ (X, T) $ is a   distal flow, then $P = \Delta = SP$.
	
	\end{remark}

\bigskip

\begin{lemma} \label{eq}
	Let  $ (X, T) $ be  a flow. Then $ SP $ is a $ T-$invariant, equivalence relation.
\end{lemma}
\begin{proof}
	Note that $SP$ is clearly reflexive, symmetric and $ T- $invariant relation. We now show the transitivity property in $SP$.
	
	Let $ (x,y), (y,z) \in SP $, and let $(a,b) \in \overline{T(x,z)} \cap \Omega$. Then there exists a $u \in J$ such that $ua = a$ and $ub = b$. Then $(a,uy,b)$ is an almost periodic point in $\overline{T(x,y,z)} \subset X \times X \times X$. But $(a,uy) \in \overline{T(x,y)}$ and $(uy,b) \in \overline{T(y,z)}$. So  $(a,uy), (uy,b) \in \Omega \cap \Delta$ and so $a=uy=b$. Thus, $(x,z) \in SP$.
\end{proof}

\begin{lemma} \label{SP equi P cond}
	Let $ (X, T) $ be  a flow. Then the following are equivalent :
	
	(i) $(x,y) \in SP $. 
	
	(ii) $ \overline{T(x, y)} \subset SP$. 
	
	(iii) $ \overline{T(x, y)} \subset P$.
\end{lemma}

\begin{proof} 
	We note that $SP$ is $T-$invariant. 
	
	For (i) $ \implies $ (ii) let $(x,y) \in SP $. Suppose $(a,b) \in  \overline{T(x, y)}$, then $ \overline{T(a, b)} \subset \overline{T(x, y)}$ and so $ \overline{T(a, b)} \cap \Omega  \subset \Delta $. Thus $\overline{T(x, y)} \subset SP \ \ \subset P$.  Then since for some minimal ideal $I$ in $E(X \times X)$ we have $p(x,y) \in \Delta$ for all $p \in I$, we must have $ \overline{T(x, y)} \cap \Omega  \subset \Delta $. Thus, $(x,y) \in SP $.
	
	(ii) $ \implies  $ (iii) is obvious.

     For (iii) implies (i) let $ (a,b)  \in  \overline{T(x,y)} \cap \Omega $. Then $ (a,b) \in \Omega \cap P $  and so $ a=b $.
     
\end{proof}

Thus, we  now can say;

\begin{theorem} \label{PEQSP}
	Let $ (X, T) $ be  a flow. Then $P$ is an equivalence relation if and only if $P = SP$.
\end{theorem}

\begin{proof}
	If $P=SP$ then by Lemma \ref{eq} $P$ is an equivalence relation.
	
	Conversely suppose $P$ is an equivalence relation and   $ (x,y) \in P $.   Then $ \overline{T(x,y)} \subset P $
	so if $ (a,b) \in \overline{T(x,y)} \cap \Omega $, then $ (a,b) \in \Omega \cap P $ and so $  a=b $. Thus  $ (x,y) \in SP $. 
\end{proof}

\bigskip

Define $SP[x] = \{y \in X : (x,y) \in SP\}$ -  the set of all points strongly proximal to $x \in X$.

\begin{lemma}
	For the minimal flow $(X,T)$ if $P$ is an equivalence relation then $SP[x] = J_Mx$, for the unique minimal ideal $M$ in $E(X)$.
\end{lemma}

\begin{proof}
	The proof follows since $E(X)$ contains a unique minimal ideal $M$, and $P[x] = J_Mx$. Since by Theorem \ref{PEQSP} $P=SP$, it follows that $P[x] = SP[x]$.
\end{proof}

\bigskip 

We can characterize $SP$ even when $P$ is not an equivalence relation.

\begin{lemma}
	Let $(X,T)$ be minimal and $x \in X$, then    $y \in SP[x]  \Longleftrightarrow$  for every minimal ideal $I$ in $E(X)$ there exists \ $  u \in J_I $, such that $y = uy =ux$,   \end{lemma}

\begin{proof}
	 Let $ y \in SP[x] $ and let $ I $ be a minimal ideal. Let $ u \in J_I $ such that $ uy=y $. Now $ (ux,uy) \in \Omega $ and so since $ (x,y) \in SP $ we have $ ux=uy=y $.
	
	Conversely suppose that for every minimal ideal $ I $ we have $ u \in J_I $ such that $ ux=uy=y $. Then if $ p \in I $ with $ px=py $  it follows that all almost periodic points in  $ \overline{T(x,y)} $ are in $ \Delta $ and so $ (x,y) \in SP $.
	\end{proof}

\begin{theorem} \label{FSP}
	Let $\pi: (X,T) \to (Y,T)$ be a factor of  flows. Then $\pi \times \pi(SP(X)) \subset SP(Y)$.
	
	Further, if $\pi$ is proximal, then $(\pi \times \pi)^{-1}(SP(Y)) = SP(X)$.
\end{theorem}
\begin{proof} Let $ (x,x') \in SP(X) $ and $ (y,y')=\pi \times \pi(x,x') $. Let $ (y_0,y_0') $ be an almost periodic point in  $ \overline{T(y,y')} $. There is an almost periodic point $ (x_0,x_0') \in  \overline{T(x,x')} $ with $ \pi \times \pi(x_0,x_0')=(y_0,y_0') $. Since $ (x,x') \in SP(X) $ we have $ x_0=x_0' $ and so $ y_0=y_0' $. Thus $ (y,y') \in SP(Y) $ and hence $\pi \times \pi(SP(X)) \subset SP(Y)$.
	
	Suppose $ \pi $ is proximal. Let $ \pi(x,x')=(y,y') \in SP(Y) $. Let $ t_i(x,x') \to (z,z') \in \Omega(X) $. Since $ (y,y') \in SP(Y) $, we have $ t_i(y,y') \to (y_0,y_0) $ and so $ \pi(z)=\pi(z') $. Thus $ (z,z') \in
	\Omega(X) \cap P(X) $ giving $ z=z' $, and  $ (x,x') \in SP(X) $.

\end{proof}

\begin{theorem} \label{R}
	For flows $(X,T)$ and $(Y,T)$, let $\pi: X \to Y$ define a proximal factor map. Then $R_\pi(X) \subset SP(X)$.
\end{theorem}
\begin{proof}
	Let $(x,x') \in R_\pi(X)$ i.e. $\pi(x) = \pi(x') $.  
	
	Let $(a,b) \in \overline{T(x, x')} \cap \Omega(X)$. Then $\pi \times \pi(a,b) \in \overline{T(\pi(x), \pi(x'))}  \subset \Delta$. Thus $(a,b) \in P(X) \cap \Omega(X)$, implying that $(a,b) \in \Delta$. 
	
	Hence $(x,x') \in SP(X)$.
\end{proof}

What about products?

\begin{theorem} \label{sp-product}
	Let $(X,T), (Y,T)$ be two flows, then $\pi_X(SP(X \times Y)) = SP(X)$ and $\pi_Y(SP(X \times Y)) = SP(Y)$ where $\pi_X$ and $\pi_Y$ are projections on $X$ and $Y$ respectively.
\end{theorem}

\begin{proof}
	If $((x,y),(x',y')) \in SP( X \times Y)$, then $ \overline{T(x, x')} \cap \Omega(X) \subset \Delta $ and $ \overline{T(y,y')} \cap \Omega(Y) \subset \Delta $. Hence $(x,x') \in SP(X)$ and $(y,y') \in SP(Y)$.
	
	Conversely let $ (x,x') \in SP(X) $, and $ (y,y') \in SP(Y) $.  We show that $ ((x,y),(x',y')) \in SP(X \times Y) $. Suppose $ t_i((x,y),(x',y')) \to ((x_0,y_0),(x_0',y_0')) \in \Omega(X \times Y) $. Then $ (x_0,x_0') \in \Omega(X) $, and $ (y_0,y_0') \in \Omega(Y)$. Also $ t_i(x,x') \to (x_0,x_0') $ so since $ (x,x') \in SP(X) $ we have $ x_0=x_0' $. Similarly $ y_0=y_0' $. Therefore $ t_i((x,y),(x',y')) \to ((x_0,y_0),(x_0,y_0)) \in \Delta$, proving that $ ((x,y),(x',y')) \in SP(X \times Y) $.
\end{proof}

\bigskip

For the flow $(X,T)$, we  recall that the pair $ (x,y) $ is said to be \emph{syndetically proximal}  if given $\epsilon > 0$  there exists a syndetic subset $ A \subset T $ such that $d(tx,ty) < \epsilon$ for all $t \in A$ as first defined in \cite{Go}.

\bigskip

\begin{theorem} \label{SP = L}
	For a flow $(X,T)$, we have $(x,y) \in SP(X) \Longleftrightarrow (x,y)$ is syndetically proximal.
\end{theorem}

\begin{proof}
	Let $(x,y) \in SP(X)$. Let $ \epsilon >0 $ and consider the $\epsilon/2-$neighbourhood $\Delta_{\epsilon/2}$ of the diagonal $\Delta$. We show that there is a syndetic $A \subset T$ such that $(ax,ay) \in \Delta_{\epsilon/2}$ for all $a \in A$.
	
	 Since $(x,y) \in SP(X)$ for every minimal ideal $M$  in $E(X)$, $M(x,y) \subset \Delta$. Also $\overline{T(x, y)} \subset SP$. For each $(a,b) \in \overline{T(x, y)}$, there exist $t_{(a,b)} \in T$ such that $t_{(a,b)}(a,b) \in \Delta_{\epsilon/2}$. Hence we have an open $U_{(a,b)} \ni (a,b)$ such that $t_{(a,b)} U_{(a,b)} \subset \Delta_{\epsilon/2}$.	 
	 
	 Now $\overline{T(x, y)}  \subset \{U_{(a,b)}: (a,b) \in \overline{T(x, y)} \}$.  Since $\overline{T(x, y)}$ is compact, we have finitely many say $(a_i,b_i) \in \overline{T(x, y)}$ for $i = 1, \ldots, n$ such that $\overline{T(x, y)} \subset \bigcup \limits_{i=1}^n \ U_{(a_i,b_i)}$,  such that for each $i$, $t_iU_{(a_i,b_i)} \in \Delta_{\epsilon/2}$ . Let $A = \{t \in T: t(x,y) \in \Delta_{\epsilon/2}\}	$, and $K = \{t_1, \ldots, t_n\}$ then we see that  $T = K^{-1}A$ and so $A$ is syndetic.

	   Let $t \in T$ then there exists $1 \leq i \leq n$ such that $t(x,y) \in U_{(a_i,b_i)}$ and so $t_i t (x,y) \in t_iU_{(a_i,b_i)} \subset \Delta_{\epsilon/2}$ and so $t_it \in A$, i.e. $t \in t_i^{-1}A$. 
	   
	   Thus	 there exists a syndetic subset $ A \subset T $ such that $d(ax,ay) < \epsilon$ for all $a \in A$ and so $(x,y)$ is syndetically proximal.
	
	Conversely suppose $ (x,y) $ is syndetically proximal. We show that $ (x,y) \in SP $.
	
	If not, there are $ t_i \in T $ with $ t_i(x,y) \to (x_0,y_0) \in \Omega(X) $ with $ x_0 \neq y_0 $.
	Then $ (x_0,y_0) \notin P(X) $. 	Let $ \epsilon >0 $ be such that $ d(tx_0,ty_0)> \epsilon $ for all $ t \in T $.
	
	Since $ (x,y) $ are syndetically proximal there is a syndetic set $ A  \subset T $ 	with $ d(ax,ay)< \epsilon/4 $ for all $ a \in A $, such that $ T=KA $ with $ K $ compact. Now $ t_i=k_ia_i $ with 	$ k_i \in K $, $ a_i \in A $, and so $ a_i(x,y) \to (z,z') $ with $ d(z,z')<\epsilon/2 $. Let $ k_i \to k \in K $ so 	$ k_ia_i(x,y) \to k(z,z') $ and  $ k_ia_i(x,y)=t_i(x,y)  \to (x_0,y_0)$. Hence $ k^{-1}(x_0,y_0)	= (z,z') $ with $ d(z,z')<\epsilon/2 $ - a contradiction.
\end{proof}

\bigskip

\begin{remark}
	The set of all $ (x,y) \in X \times X $ that are syndetically proximal forms an invariant, equivalence relation $L$ which was first studied in \cite{C}. In view of Theorem \ref{SP = L}, our Theorem \ref{SP equi P cond} happens to be same as Lemma 4 and Lemma 5 in \cite{C}. Thus $SP = L$.
\end{remark}

\bigskip

We recall that the flow  $(X,T)$ is said to be \emph{locally almost periodic at $ x  \in X$}  if, whenever $ U $ is a neighbourhood of $ x $, there is a neighbourhood $ V $ of $ x $ and a syndetic subset $ A \subset T $ such that
$ AV \subset U $. We say that $(X,T)$ is \emph{locally almost periodic (lap)} if it is locally almost periodic at every point.

 Let  $(X,T)$ and $(Y,T)$ be minimal flows. Recall that $ \pi:X \to Y $ is \emph{highly proximal (HP) } if $ x \in \pi^{-1}(y) $ and $ p   $ in any minimal ideal $I \subset E(X)$, then $ p \circ \pi^{-1}(y)=\{px\} $, i.e., if for some $ y \in Y  $ there is a net $ \{t_n\} $ in $ T $ and a point    $ x \in X $ such that   $ t_n \pi^{-1}(y) \to \{x\} $ in the Hausdorff topology   of $ 2^X $. 
 
  It is known that in such a case  $ (X,T) $ is locally almost periodic if and only if $ (Y,T) $ is locally almost periodic. We refer to \cite{AG, AM, GH} for all related theory and state an  important theorem in this direction,

\begin{theorem} \cite{AM} \label{lap+hp}
	A minimal flow is locally almost periodic if and only if it is a highly proximal extension of an equicontinuous flow.
\end{theorem}

which has an interesting corollary,

\begin{corollary} \cite{AM} \label{lap=1-1}
	Let $ (X, T) $  be a minimal flow. Then $ (X, T) $ is locally almost periodic if and only if
	there is an equicontinuous minimal flow $ (Z, T) $, a homomorphism $ \pi : X  \to Z $, and a residual
	subset $ E $ of $ Z $ such that for all $ z \in E $, $ \pi^{-1}(z) $ consists of a single point.
\end{corollary}

And we observe that:

\begin{theorem}
	For the flow $(X,T)$ let all minimal subflows be locally almost periodic flows. Then $P(X) = SP(X)$. 
\end{theorem}

\begin{proof}
	Let $(x,y) \in P(X)$, and $z \in X$ be almost periodic with $p \in E(X)$  such that $px = z = py$. Then $(z,z) \in  \Omega(X)$ and is locally almost periodic.
	
	Let $\epsilon > 0$ and $W$ be an $\epsilon/2-$neighbourhood of $(z,z) \in X \times X$. Then there exists an open $V \ni (z,z)$ and a syndetic $A \subset T$ such that $AV \subset W$.
	
	Since $(z,z) \in \overline{T(x,y)}$, there exists $t \in T$ such that $t(x,y) \in V$. 
	
	Now $At(x,y) \subset AV \subset W$ with $At$ syndetic, and so  $d(ax,ay) < \epsilon, \ a \in At$. Thus, $(x,y)$ is sydetically proximal and so $(x,y) \in SP(X)$.
	
\end{proof}

We note that it is known that for minimal locally periodic flows, $SP(X) = P(X)$, i.e., $P(X)$ is an equivalence relation. However, our theorem holds for non minimal flows. We now see that there are many such non minimal flows for which $P(X)$ is an equivalence relation.

\begin{example}
	Let $ (X_i,T) $  be locally almost periodic (lap) minimal flows, for $1 \leq i \leq n$, $n \in \N$. Then there are equicontinuous minimal flows $ (Z_i,T) $, $1 \leq i \leq n$, such that $ (X_i,T) $ are HP extensions of $ (Z_i,T) $  respectively, where $ \pi_i:X_i \to Z_i $, $1 \leq i \leq n$, are the homomorphisms. Let $ x_i \in X_i $ be such that $ (x_1, \ldots, x_n) $ is an almost periodic point of $\prod_{i=1}^{n} X_i$ and let $ K= \overline{T(x_1, \ldots, x_n)} $ so $ K $ is minimal. We invoke Theorem \ref{lap+hp} to show that $ (K,T) $ is lap. Let  $z_i \in Z_i$ be such that $\pi(x_i) = z_i$, $1 \leq i \leq n$ and let $L = \overline{T(z_1, \ldots, z_n)} $. Now $\prod_{i=1}^{n} \pi_i:K \to L$, and $L$ is a minimal subset of $\prod_{i=1}^{n} Z_i $, so $ L $ is equicontinuous. Suppose $ t_i \to p $  in $ \cM $. Then $t_i \prod_{i=1}^{n} p_i^{-1} (z_1, \ldots, z_n) \to \{(px_1, \ldots, px_n)\}=\{p(x_1, \ldots, x_n)\}$ and hence $ \prod_{i=1}^{n} \pi_i $ is HP giving $  K $ to be lap, and $P(\prod_{i=1}^{n} X_i)$ is an equivalence relation. 
\end{example}

\begin{definition}
The \emph{weakly distal  relation} is defined as 

$$WD(X) = WD= \{(x,y) \in X \times X: \overline{T(x, y)} \cap \Omega \nsubseteq \Delta \}. $$

\end{definition}

Clearly, $D \subseteq WD$.

Define $WD[x] = \{y \in X : (x,y) \in WD\}$ -  the set of all points weakly distal  to $x \in X$.

\bigskip

We have $(x,y) \in SP \ \Longleftrightarrow \ (x,y) \notin WD \ \text{and} \ (x,y) \in WD \ \Longleftrightarrow \ (x,y) \notin SP$. And $X \times X = SP \cup WD$ is a disjoint union.

$WD$ is  symmetric, $ T- $ invariant
relation, and cannot be reflexive since $WD \cap \Delta = \emptyset$.

\bigskip

We see that $WD = (P \setminus SP) \cup D$. So  $P = SP  \Longleftrightarrow WD = D$.

\begin{theorem} \label{proxiFWD}
	For flows $(X,T)$ and $(Y,T)$ let $\pi: X \to Y$ be a  factor map. Then $({\pi \times \pi})^{-1}(WD(Y)) \subset WD(X)$.
	
	Furthermore, if $\pi$ is proximal then $\pi \times \pi(WD(X)) \subset WD(Y)$.
\end{theorem}

Note that this follows by complementation of Theorem \ref{FSP}.

\begin{theorem} \label{wd-product}
	Let $(X,T), (Y,T)$ be two flows and $\pi_X$ and $\pi_Y$ be projections of $X \times Y$ on $X$ and $Y$ respectively. Then either one of the following or both hold true:
	
	1. $\pi_X(WD(X \times Y)) \subset WD(X)$ 
	
	2. $\pi_Y(WD(X \times Y)) \subset WD(Y)$,

\end{theorem}

\begin{proof}
	If $((x,y),(x',y')) \in WD( X \times Y)$, then we have either $ \overline{T(x, x')} \cap \Omega(X) \nsubseteq \Delta $ or $ \overline{T(y,y')} \cap \Omega(Y) \nsubseteq \Delta $ or both. Hence $(x,x') \in WD(X)$ or $(y,y') \in WD(Y)$ or both hold true.

\end{proof}

\begin{lemma}
	For the minimal flow $(X,T)$ if $P(X)$ is an equivalence relation then $WD[x] = X \setminus J_Mx$, for the unique minimal ideal $M$ in $E(X)$.
\end{lemma}

\begin{proof}  From Theorem \ref{PEQSP},  $P = SP$. The proof follows as $WD = X \times X \setminus SP$.
\end{proof}

\begin{remark}
	Let $(Z,T)$ be a subsystem of $(X,T)$, then $SP(Z) \subseteq SP(X)$ and $WD(Z) \subseteq WD(X)$.
\end{remark}

\bigskip

\begin{theorem} \label{prodSP} Let $(\{(X_\alpha, T): \alpha \in \Lambda\})$ be a family of flows with action by the same group $T$, and let $X = \prod \limits_{\alpha \in \Lambda} X_\alpha$.
	
	Then	$(x,y) \in SP(X) \ \Longleftrightarrow (x_\alpha, y_\alpha) \in SP(X_\alpha) $ for all $\alpha \in \Lambda$.
\end{theorem}
\begin{proof}
	Note that $(x,y) \in SP(X) $ implies that $I_\nu(x,y) \subset \Delta$, for every minimal ideal $I_\nu$ in $E(X)$. But this means that $I_{\alpha,\nu}(x_\alpha,y_\alpha) \subset \Delta$, for every minimal ideal $I_{\alpha,\nu}$ in $E(X_\alpha)$ for all $\alpha \in \Lambda$. Thus $(x_\alpha, y_\alpha) \in SP(X_\alpha) $ for all $\alpha \in \Lambda$.
\end{proof}

\begin{theorem} 
	Let $(\{(X_\alpha, T): \alpha \in \Lambda\})$ be a family of flows with action by the same group $T$, and let $X = \prod \limits_{\alpha \in \Lambda} X_\alpha$. 
	
	Then	$ (x,y) \in WD(X) $  if and
	only if $(x_\alpha, y_\alpha) \in WD(X_\alpha) $ for some $\alpha$.
\end{theorem}
\begin{proof}
	We see that $(x,y) \in WD(X) $ implies that $\overline{T(x,y)} \cap \Omega(X) \nsubseteq \Delta $ which gives  $\overline{T(x_\alpha,y_\alpha)} \cap \Omega(X_\alpha) \nsubseteq \Delta $ for some $\alpha$. Conversely, if  $\overline{T(x_\alpha,y_\alpha)} \cap \Omega(X_\alpha) \nsubseteq \Delta $ for some $\alpha$ then $\overline{T(x,y)} \cap \Omega(X) \nsubseteq \Delta $. 
\end{proof}

\bigskip

\bigskip

\bigskip

\section{Proximal and Strongly Proximal Sets}

\begin{definition}
	Let $ (X, T) $ be a 	flow. A  set $A   \subset X$ is called  a \emph{proximal set} if there exists $p \in E(X)$ such that $pA$ is a singleton set, i.e. there exists $z \in X$ such that $pA = \{z\}$.
	
\end{definition}	

\begin{theorem}	\label{like} Let $ (X, T) $ be a 	flow and $A \subset X$. Then the following are equivalent:
	
	1. 	$A$ is a proximal set.
	
	2. Every    $\bar{x} \in X^{|A|}$ with $range(\bar{x}) = A$   is proximal in $X^{|A|}$, where  $ |A| $ is the cardinality of $ A $, i.e. $p\bar{x} \in \Delta \subset  X^{|A|}$ for some $p \in E(X)$.
	
	3.  If we write $A = \ \{x_\alpha: \alpha \in \Lambda\} $ with $|\Lambda| = |A|$, then there is a  $ p \in E(X) $ such that $ p((x_\alpha)_{\alpha \in \Lambda}) = \underbrace{(z, \ldots, z)}_{|\Lambda| \ \text{times}}  $ for some $z \in X$. 
	
\end{theorem}

\begin{proof}
	For 1. $ \implies $ 2.	let $A$ be a proximal set, then there exists $p \in E(X)$ such that $pA = \{z\}$ is a singleton. Now consider the   tuple $\bar{x} \in X^{|A|}$ with $range(\bar{x}) = A$. Then $p \bar{x} = \underbrace{(z, \ldots, z)}_{|A| \ \text{times}}  $. Thus $p\bar{x} \in \Delta \subset  X^{|A|}$ for some $p \in E(X)$. 
	
	For 2. $ \implies $ 3. if $A = \ \{x_\alpha: \alpha \in \Lambda\} $ with $|\Lambda| = |A|$, there is a  $ p \in E(X) $ such that $ p((x_\alpha)_{\alpha \in \Lambda}) = \underbrace{(z, \ldots, z)}_{|\Lambda| \ \text{times}}   = p\bar{x}  \in X^{|A|}$. 
	
	For 3. $ \implies $ 1. we note that for  $A = \ \{x_\alpha: \alpha \in \Lambda\} $ with $|\Lambda| = |A|$, if there is a  $ p \in E(X) $ such that $ p((x_\alpha)_{\alpha \in \Lambda}) = \underbrace{(z, \ldots, z)}_{|\Lambda| \ \text{times}}  $ for some $z \in X$. Then  we have $pA = \{z\}$ - a singleton. Thus $A$ is a proximal set.
\end{proof}

Each statement in the theorem above can be taken to be the definition of a proximal set. 

\bigskip

\begin{remark}
	Note that a subset of a proximal set is proximal. 
	
In particular if $ (X,T) $ is a proximal flow, i.e. $P = X \times X$, then $X$ itself is a proximal set, as the minimal idempotents in E(X) are constants and so for any $u \in J$, the set $uX$ will be a singleton. 	
\end{remark}

\begin{example}
	Let $(X,T)$ be a flow such that $P(X)$ is an equivalence relation. For $x \in X$, if $y \in P[x]$ then $P[x] = P[y]$. Thus for $x,y \in X$ either $P[x] = P[y]$ or $P[x] \cap P[y] = \emptyset$. Thus each $P[x]$ are proximal sets.
	
	When $P(X)$ is not an equivalence relation, $X$ still admits  proximal sets. See Example \ref{CC} where $A$ and $B$ are two distinct  proximal sets.
\end{example}

\begin{remark} Note that every finite subset of a proximal set is proximal. However, if we assume that every finite subset of a set is proximal than the set need not be proximal. We consider the $\emph{2}-$shift $(\emph{2}^{\Z},\sigma)$. We consider the metric $d$ on $\emph{2}^{\Z}$ defined as $d(x,y) = \frac{1}{2^{n+1}}$ whenever the blocks $x_{[-n,n]} = y_{[-n,n]}$.
	
	We note that for $ \epsilon > 0 $, there exists $ t \in \N $ such that $ d(x,\mathbf{0}) < \epsilon $, where $ \mathbf{0} = 0^\infty.0^\infty$ - the zero bi-infinte sequence, if $ x_i = 0 $ for $ -t \leq i \leq t $.
	
	\bigskip
	
	Now let $ A = \{x^M \in \emph{2}^{\Z}: M \in \N\} $ where $ x^M = 0^\infty.1^M 0^\infty. $
	
	\bigskip
	
	Then for any finite subset $ F = \{x^{K_1}, \ldots, x^{K_l}\} \subset A $, let $x^K \in F$ be such that $K = \max \{K_1, \ldots, K_l\}$. Then  we have $ d(\sigma ^{-t-K}x, \mathbf{0}) < \epsilon $ for all $x \in F$. Since $ \sigma^{-n}F $ continues to be in  a neighbourhood of $ \mathbf{0} $ for $ -n < -t -K $ we have a $  p \in E(\emph{2}^{\Z}) $ such that $ \sigma ^{-n} \to p $ and $ pF = \{\mathbf{0}\} $. And so $ F $ is a proximal set.
	
	[We use the fact that if $ pF = \{\mathbf{0}\} $ then for $ \sigma^r \to p $ we should have some $ \sigma ^sF $ in the neighbourhood of $ \mathbf{0} $.]
	
	However, $ A $ is not proximal as there exists no $k \in \N$ for which $ d(\sigma ^{-k}x, \mathbf{0}) < \epsilon $ for all $x \in A$. In fact for any $n \in \N$, there is a $N \in \N$ such that $ d(\sigma ^{N}x^i, \mathbf{0}) < \epsilon $ for all $i < n$, but $ d(\sigma ^{N}x^j, \mathbf{0}) > \epsilon $ for all $j > > n$ for elements in $A$.
	
	In this sense our definition of proximal set differs from that of \cite{AUS} given on page 92.

\end{remark}

We note this trivial result.

\begin{lemma}
	Let $(X,T)$ be a flow. Then for every proximal set $A \subset X$, we have $tA$ is also proximal for all $t \in T$. 
\end{lemma}

Suppose $ P $ is an equivalence relation and $ A $ is a proximal set. Then note that  $ A \subset P[x] $ for all $ x \in A $.

\begin{theorem}
	For the flow $(X,T)$ let $P$ be an equivalence relation $\Longleftrightarrow$ $rA$  is a proximal set for each proximal set $A$ and for $r \in E(X)$. 
\end{theorem}
\begin{proof}
 Note that a proximal set	$ A  \subset P[x] $, for some $x \in A$,  since $ P $ is an equivalence relation. And $ rA \subset P[rx]$ for all  $r \in E(X)$.
 
 Conversely, let $ (x,y) \in P $. Then $ A=\{x,y\} $ is a proximal set. If $  r \in E $, then $ rA=\{rx,ry\} $ is also a proximal set. It follows that $ \overline{T(x,y)} \subset P $, so $ P $ is an equivalence relation.
 
\end{proof}

\begin{theorem}
	Let $ (X,T) $ be a flow, and let $A \subset X$ be a proximal set. Then  there is a minimal  ideal $ I $ in $ E(X) $ such that $pA $ is a singleton set	for all $ p \in I $.
\end{theorem}

\begin{proof}
	Let $A \subset X$ be proximal. Then  there exists $ p \in E(X) $  such that $ pA = \{z\}  $ for some $z \in X$.
	
	Let $\mathcal{C} = \{p \in E(X): pA \ \text{is a singleton set} \}$. Then we note that $\mathcal{C} \neq \emptyset$ and $\mathcal{C} \subset E(X)$ is a left ideal. Hence it must contain a minimal (left) ideal $I$. So  $pA$ is a singleton set, for all $p \in I$. 
\end{proof}

\begin{theorem}
	For flows $(X,T)$ and $(Y,T)$ let $\pi: X \to Y$ define a  proximal factor map. Then  $\pi^{-1}(y)$ is a proximal set, for all $y \in Y$.
\end{theorem}

\begin{proof} Let $x,x' \in \pi^{-1}(y)$ and  $ I $ be a minimal  ideal then for $ p \in I $, we have 	$ (px,px') \in \Omega(X)$. But $ \pi(px)=\pi(px') $  and so $ (px,px') \in P(X) $ and thus $ px=px' = z $ (say). Hence $ px=z $ for all $ x \in pi^{-1}(y) $ making it a proximal set.	 
\end{proof}

\begin{definition}
	Let $I \subset E(X)$ be a minimal ideal, then $A_I \subset X$ is called a \emph{$I-$proximal set}, if   we have $pA_I$ is a singleton for all $p \in I$ and $A_I$ is maximal with respect to this property.
\end{definition}

We note that there can be more than one $I-$proximal sets for a given minimal ideal $I$ in $E(X)$.

\begin{example}
	We refer to Example \ref{MT}. Here there are two minimal ideals $I$ and $K$. There are infinitely many distinct $I-$proximal sets $\{a,b\}$,  $\{\bar{a},\bar{b}\}$ and their orbits, and infinitely many $K-$proximal sets $\{a,\bar{b}\}$, $\{\bar{a},b\}$ and their orbits.
	
	Let $(X,T)$ be a flow such that $P(X)$ is an equivalence relation. For $x \in X$, and the only minimal ideal $I$ in $E(X)$, each $P[x]$ are $I-$proximal sets as can be seen in Example \ref{F}.
	
\end{example}

\begin{theorem}
	For the flow $(X,T)$, let $I$ be a minimal ideal in $E(X)$ and $A_I$ be a $I-$proximal set. Then $A_I$ contains an almost periodic point.
	
	Further, for every $u \in J_I$, $ua \in A_I$ for all $a \in A_I$.
\end{theorem}

\begin{proof}
	Suppose that $A_I$ does not contain any almost periodic point. Then for all $a \in A_I$ and all $u \in J_I$, we have $ua \notin A_I$. But then if $pA_I = \{z_p\}$ then $p(ua) = pa = z_p$ for every $p \in I$, $u \in J_I$ and $a \in A_I$. 	Hence, $A_I \cup \{ua\}$ is a proper superset of the $I-$proximal set $A_I$ contradicting the maximality of $ A_I $.  Thus $ua \in A_I$. 
	
	Hence for every $u \in J_I$ and every $a \in A_I$, we have $ua \in A_I$. 
\end{proof}

\begin{remark}
	Note that for any $a,b \in A_I$ and $u \in J_I$, we have $ua = ub = z_u$. And so $z_u \in A_I$. Also $pz_u = pua = pa = z_p$ for all $p \in I$. Thus $IA_I = I\{z_u\}$ is a minimal subset of $X$, with $uz_u = z_u$.
	
	Also when $P(X) = X \times X$, then for the only minimal ideal $I$ in $E(X)$, we have $X$ is the only $I-$proximal set. 
	
\end{remark}

In general, a $I-$proximal set need not be closed.

\bigskip

Fix  a  minimal ideal $ I $ in $ E(X)$.

\begin{theorem}
	Let $(X,T)$ be a flow with $I$ a minimal ideal. Then $pA_I \neq pB_I$ for all $p \in I$ and $I-$proximal sets $A_I \neq B_I$.
\end{theorem}
\begin{proof}
	Let  $A_I \neq B_I$ be $I-$proximal. Then there exists $q \in I$ with $qA_I \neq qB_I$. Note that otherwise $A_I \cup B_I$ is $I-$proximal contradicting the maximality of both $A_I$ and $B_I$. 
	
	Suppose for some $p \in I$ we have $pA_I = pB_I$. Then $upA_I = upB_I$ for $u \in J_I$. Since $uI$ is a group, there exists $up' \in uI$ such that $up'(up) = u -$ the identity in $uI$. Thus, $uA_I = uB_I$ and so $pA_I = puA_I = puB_I = pB_I$ for every $p \in I$. This gives a contradiction.
	
	Thus $pA_I \neq pB_I$ for all $p \in I$. 
	
\end{proof}

\begin{theorem}
	Let $(X,T)$ be a flow. Then for every $x \in X$, there exists an $I-$proximal set $A_I$ such that $x \in A_I$.
\end{theorem}

\begin{proof}
	For $x \in X$ consider the set $B = \{x\} \cup \{ux: u \in J_I\}$. Then $pB =\{px\}$ for all $p \in I$ and so is a proximal set. Let $ \mathcal{A} =\{ A \subset X: B \subset A$ and $pA $ is a singleton for each $p \in I\}$.
	
	Consider a totally ordered set $\{A_i\}$ of elements in  $\mathcal{A}$, then $ \cup_i A_i \in \mathcal{A}$. Thus every chain in $\mathcal{A}$ is bounded above and so by Zorn's lemma $\mathcal{A}$ has a maximal element (say) $A_I$. Thus $A_I$ is an $I-$proximal set with $x \in A_I$.
\end{proof}

 Define a relation `$\simeq_I$' on $X$ as $x \simeq_I y$ if and only if there exists an $I-$proximal set $A_I$  such that $x,y \in A_I$.

\begin{lemma}
	For the flow $(X,T)$, the relation $\simeq_I$ is an equivalence relation.
\end{lemma}

\begin{proof} We only need to look into the transitivity of the relation $\simeq_I$.
	
	Let $x \simeq_I y$ and $y \simeq_I z$. Then there exists $I-$proximal sets $A_I$ and $B_I$ such that $x,y \in A_I$ and $y,z \in B_I$. Then $px = py$ and $py = pz$ for all $p \in I$. This gives $px = pz$ for all $p \in I$ implying that $x,y,z \in A_I$ or $B_I$ by the maximality of $A_I$ and $B_I$. Thus $x \simeq_I z$ and so $\simeq_I$ is an equivalence relation.
\end{proof}

\begin{remark}
	For the flow $(X,T)$ let $P(X)$ be closed and hence an equivalence relation. Then  for the only minimal ideal $I$ in $E(X)$, we note that $\simeq_I$ is same as the relation $P$ and hence an icer.
	
	Again if $I$ and $K$ are two distinct ideals in $E(X)$ then $\simeq_I$ need not be the same as $\simeq_K$.
\end{remark}

\begin{definition}
	Let $ (X, T) $ be a 	flow. A  set $A   \subset X$ is called  a \emph{strongly proximal set} if for every minimal ideal $M$ in $E(X)$ and every $p \in M$ the set $pA$ is a singleton set, i.e. there exists $z \in X$ such that $pA = \{z\}$.
	
\end{definition}

Analogous to Theorem \ref{like} we have:	

\begin{theorem}	Let $ (X, T) $ be a 	flow and $A \subset X$. Then the following are equivalent:
	
	1. 	$A$ is a strongly proximal set.
	
	2. For every    $\bar{x} \in X^{|A|}$ with $range(\bar{x}) = A$, where  $ |A| $ is the cardinality of $ A $, and every $p \in M$ for every minimal ideal $M$ in $E(X)$ the tuple $p\bar{x} \in \Delta \subset  X^{|A|}$.
	
	3.  If we write $A = \ \{x_\alpha: \alpha \in \Lambda\} $ with $|\Lambda| = |A|$, for  every $p \in M$ for every minimal ideal $M$ in $E(X)$, we have $ p((x_\alpha)_{\alpha \in \Lambda}) = \underbrace{(z, \ldots, z)}_{|\Lambda| \ \text{times}}  $ for some $z \in X$. 
	
\end{theorem}

\begin{proof}
	For 1. $ \implies $ 2.	let $A$ be a strongly proximal set, then and every $p \in M$ for every minimal ideal $M$ in $E(X)$ the set $pA = \{z\}$ is a singleton. Now consider the   tuple $\bar{x} \in X^{|A|}$ with $range(\bar{x}) = A$. Then $p \bar{x} = \underbrace{(z, \ldots, z)}_{|A| \ \text{times}}  $. Thus $p\bar{x} \in \Delta \subset  X^{|A|}$ for and every $p \in M$ for every minimal ideal $M$ in $E(X)$.
	
	For 2. $ \implies $ 3. if $A = \ \{x_\alpha: \alpha \in \Lambda\} $ with $|\Lambda| = |A|$, and every $p \in M$ for every minimal ideal $M$ in $E(X)$ $ p((x_\alpha)_{\alpha \in \Lambda}) \to \underbrace{(z, \ldots, z)}_{|\Lambda| \ \text{times}}   = p\bar{x}  \in X^{|A|}$. 
	
	For 3. $ \implies $ 1. we note that for  $A = \ \{x_\alpha: \alpha \in \Lambda\} $ with $|\Lambda| = |A|$, and every $p \in M$ for every minimal ideal $M$ in $E(X)$ the tuple $ p((x_\alpha)_{\alpha \in \Lambda}) = \underbrace{(z, \ldots, z)}_{|\Lambda| \ \text{times}}  $ for some $z \in X$. Thus  we have $pA = \{z\}$ - a singleton, i.e. $A$ is a strongly proximal set.
\end{proof}

Each statement in the theorem above can be taken to be the definition of a strongly proximal set. 

\bigskip

\begin{remark}
	Note that a subset of a strongly proximal set is strongly proximal. 
	
	In particular if $ (X,T) $ is a proximal flow, i.e. $P = X \times X(= SP)$, then $X$ itself is a  strongly proximal set. 
\end{remark}

Consider the set $\mathfrak{S}$ of all strongly proximal sets. Then every chain in $\mathfrak{S}$ has an upper bound and so has a maximal. These maximal elements are called \emph{maximal strongly proximal sets}.

Note that if $A$ and $B$ are maximal strongly proximal sets then $A \cap B = \emptyset$. Else $A \cup B$ becomes strongly proximal.

\begin{theorem} 
	For the collection $\mathfrak{M}$ of minimal ideals in $E(X)$, let $\{A_M : M \in \mathfrak{M}\}$ be a collection of $M-$proximal sets in $X$.	Then $\bigcap \limits_{M \in \mathfrak{M}} A_M$ is a maximal strongly proximal set.
\end{theorem}
\begin{proof}
	Let $ A = \bigcap \limits_{M \in \mathfrak{M}} A_M$. Then for every minimal ideal $M \in \mathfrak{M}$ and every $p \in M$ the set $pA$ is a singleton set, and hence $A$ is a strongly proximal set. We now show that it is maximal. Suppose $A \cup \{x\}$ is also strongly proximal for $x \notin A$. Then for every minimal ideal $M \in \mathcal{M}$ and every $p \in M$ the set $p(A \cup \{x\})$ is a singleton set. But this means that for every $M \in  \mathcal{M}$ the set $A_M \cup \{x\}$ should be a $M-$proximal set, contradicting the maximality of $A_M$.
	
	Hence $A$ is a maximal strongly proximal set.
\end{proof}

\begin{remark}
	The above theorem implies that for every $I \in \mathfrak{M}$, every $I-$proximal set contains a maximal strongly proximal set.
\end{remark}

\begin{example}
	We refer to Example \ref{MT}. Here there are two minimal ideals $I$ and $K$ and we have minimal idempotents $u_1, u_2 \in I$ and $V_1, v_2 \in K$. There we consider  two $I-$proximal sets $\{a,b\}$ and $\{\bar{a},\bar{b}\}$, and two $K-$proximal sets $\{a,\bar{b}\}$ and $\{\bar{a},b\}$. So now there are two choices for $A_I$ and $A_K$ and this gives four possibilities for $A_I \bigcap  A_K$ - namely $\{a\}, \ \{b\}, \ \{\bar{a}\}, \ \{\bar{b}\}$, which are maximal strongly proximal sets.
\end{example}

 Define a relation `$\equiv$' on $X$ as $x \equiv y$ if and only if there exists a maximal strongly proximal set $A$  such that $x,y \in A$.

\begin{lemma}
	For the flow $(X,T)$, the relation $\equiv$ is an equivalence relation.
\end{lemma}

\begin{proof} We only need to look into the transitivity of the relation $\equiv$.
	
	Let $x \equiv y$ and $y \equiv z$. Then there exists maximal strongly proximal sets $A$ and $B$ such that $x,y \in A$ and $y,z \in B$. Then $A \cap B \neq \emptyset$ and by  the maximality of $A$ and $B$, we must have $ A = B$. Thus $x \equiv z$ and so $\equiv$ is an equivalence relation.
\end{proof}

\begin{theorem}
	For the flow $(X,T)$, let  $A$ be a maximal strongly proximal set. Then for all $u \in J$ we have the singleton set $uA \subset A$.
\end{theorem}

\begin{proof} Let $I$ be a minimal ideal in $E(X)$. Suppose that $uA = \{a_u\}$ (say) for some $u \in J_I$. Then for all $p \in I$, we have $pA = p\{a_u\} = \{a_p\}$ (say). We note that $ua_u = a_u$ and so $p(A \cup \{a_u\}) = \{a_p\}$ and so by maximality of $A$ we must have $a_u \in A$. 
	
	Since $I$ was arbitrary we have $a_u \in A$ for all  all $u \in J$.
	 
\end{proof}

\bigskip

\section{Weakly Distal Flows}

We say that the relation $WD$ is full if $WD = X \times X \setminus \Delta$, i.e. $SP = \Delta$. 

\begin{definition}
	The  flow $ (X,T) $ is called \emph{weakly distal} $\Longleftrightarrow$ the relation $WD$ is full \ $\Longleftrightarrow \ SP = \Delta$.
\end{definition} 

We note some analogy here. If $P = \Delta$ the flow is distal and if $SP = \Delta$ then the flow is weakly distal. In general since $SP \subset P$, we can say that all distal flows are also weakly distal.

\begin{remark}
	The Examples \ref{CC}, \ref{MT}, \ref{PL} and \ref{Chacon} described in the next section have $SP = \Delta$ i.e. a full $WD$ and so are weakly distal flows. Whereas Examples  \ref{F} and \ref{Shift} have $SP \neq \Delta$, i.e. $WD$ is not full here and so are not weakly distal.

\end{remark}

\begin{lemma}
	A subflow of a weakly distal flow is weakly distal.
\end{lemma}

\begin{proof} The proof follows since $SP(Z) \subseteq SP(X)$ whenever  $(Z,T)$ is a subflow of $(X,T)$. \end{proof}

\bigskip

What happens with factors and extensions? 

\begin{lemma}
	For a flow $(X,T)$, if $SP(X) (\neq \Delta)$ is closed then the canonical factor $(X/SP(X), T)$ is weakly distal.
\end{lemma}

\begin{proof}
	For a closed $ SP(X) $ let $ Y= X/SP(X) $ and suppose $ (y,y') \in SP(Y) $. We show that  $ y=y' $. Let $\pi: (X,T) \to (Y,T)$ be  the natural canonical factor. 
	
	Let $ \pi(x,x')=(y,y') $. Let $ p(x, x') = (x_0,x_0') \in \Omega(X) $  for some $p \in E(X)$. Let $\theta: E(X) \to E(Y)$ be the natural homomorphism with $ \theta(p) \in E(Y) $. Let $ \theta(p)(y,y') =  (y_0,y_0') \in \Omega(Y) $.
	Since $ (y,y') \in SP(Y) $ we have $ y_0=y_0' $. So $ \pi(x_0)= \pi(x_0') $. Hence $ (x_0,x_0') \in SP(X) \subset P(X) $.
	Thus $ (x_0,x_0') \in P(X) \cap \Omega(X) $ and so $ x_0=x_0' $. But this means that $ (x,x') \in SP(X)$ and so $ 	\pi(x)=\pi(x') $ which means $ y=y' $.
\end{proof}

Note that  for a flow $(X,T)$ if $P(X)$ is closed, then $P = SP$ and so $D = WD$ and so the canonical factor $(X/P(X), T)$ is distal.

\begin{remark} Note that if $\pi: (X,T) \to (Y,T)$ is  proximal then, as pointed out in Theorem \ref{R}, $ \pi(x)=\pi(x') $ for $x \neq x'$	then $ (x,x') \in SP(X) $. So $\pi$ will be an isomorphism if $(X,T)$ is weakly distal. Though, a  weakly distal flow need not be preserved under proximal  extensions. 
\end{remark}

\begin{definition}
	For flows $(X,T)$ and $(Y,T)$,  the extension $\pi: X \to Y$  is called  \emph{weakly distal } if $\pi(x) = \pi(x')$ then $(x,x') \in WD(X)$.
\end{definition}

\begin{example} \label{adic}
	We note the following example from \cite{AUS} (page 241):
	
	Consider the endomorphism $ H $ on  $ \{ 0 , 1 \}^\Z  $ defined
	by $ Hw(n) = w(n)+w(n+1) $ for all $w =(w(n))_{n \in \Z}$. Then $H$ maps the Morse-Thue minimal flow $(X, \sigma)$ (Example \ref{MT}) onto a minimal	flow $ (X_H, \sigma) $.
	
	Now $ H(x_l) = H(x_2) $ if and only  if $ x_2 = x_1^* $  (where $w^*$ denotes the bi-infinte sequence $w$ with $0$ and $1$ interchanged).
	This $ H $ is an	equicontinuous extension. And $ X_H $ is an almost one-to-one (hence proximal) extension of the ``adding machine" -- the equicontinuous	flow $ ( D , Add ) $ on the dyadics.
	
	Since $H: X \to X_H$ and $L: X_H \to D$ are factor maps, this gives a distal flow $(D, Add)$ as a factor of the weakly distal flow $(X, \sigma)$. Hence the extension $L \circ H$ is weakly distal.
\end{example}

\begin{lemma}
	A weakly distal extension of a weakly distal flow is weakly distal.
\end{lemma}

\begin{proof}
	For flows $(X,T)$ and $(Y,T)$, let $\pi: X \to Y$ be weakly distal. Let $ (x,x') \in SP(X) $. Then $ (\pi(x),\pi(x')) \in SP(Y) $.
	But since $ (Y,T) $ is weakly distal, we have $ \pi(x)=\pi(x') $. But since $ \pi $ is weakly distal, we must have $ x=x' $. Thus $SP(X) =\Delta$ and $(X,T)$ is weakly distal.
\end{proof}

	\bigskip

\begin{remark}
	
	Note that  in the Example \ref{adic}, the Morse minimal flow $X$ is an equicontinuous extension of $X_H$ which is an almost one-to-one extension of an equicontinuous flow  $D$. Clearly $X_H$ is not weakly distal as $SP(X_H) \neq \Delta$. Thus,  weakly distal flows need not have   a  weakly distal factor.

Again for the Example \ref{Shift}, we note that though $(X,G)$ is not weakly distal since $SP$ is closed,  so the canonical factor $(X/SP,G)$ is weakly distal.

This also illustrates that in general, extensions of weakly distal flows need not be weakly distal.\end{remark}

\bigskip

\begin{theorem}  Let $(\{(X_\alpha, T): \alpha \in \Lambda\})$ be a family of weakly distal flows with action by the same group $T$, and let $X = \prod \limits_{\alpha \in \Lambda} X_\alpha$.
	
	Then $(X,T)$ is weakly distal. \end{theorem}
\begin{proof}  It follows from Theorem \ref{prodSP} that $(x,y) \in SP(X) \ \Longleftrightarrow (x_\alpha, y_\alpha) \in SP(X_\alpha) $ for all $\alpha \in \Lambda$. Since  $SP(X_\alpha) = \Delta $ for all $\alpha \in \Lambda$, we have $SP(X) = \Delta$ and so $(X  T)$ is weakly distal. \end{proof}

\begin{corollary}
	Let $(X,T)$ be a weakly distal flow. Then $(X^\Lambda, T)$ is also weakly distal for any cardinality $\Lambda$.
\end{corollary}

We look into the enveloping semigroups of weakly distal flows.

\begin{theorem}
	For the weakly distal flow $(X,T)$, its enveloping semigroup flow $(E(X), T)$ is also weakly distal.
\end{theorem}
\begin{proof}
	We note that the flow $(X^X,T)$ is weakly distal and since $E(X)$ being a subflow of $X^X$ is also weakly distal. 
\end{proof}

\bigskip

We note that all distal flows are also weakly distal, and thus weakly distal flows are a bigger class of flows than distal flows. However on the basis of the enveloping semigroups we can distinguish betwen these classes of flows.

\bigskip

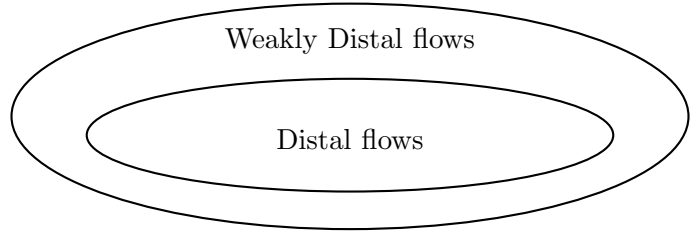
\begin{center}
	\begin{tikzpicture}
		\draw [thick] (0,-0.5) arc (-90:270:4.5cm and 1.5cm);
		\draw [thick] (0,0) arc (-90:270:3.5cm and 0.75cm);
		\node [yshift=.7cm] (0,0) {$\text{Distal} \ \text{flows}$};
		\node [yshift=2cm] (0,0) {$ \text{Weakly Distal} \ \text{flows}$};
		\end{tikzpicture}
\end{center}

\bigskip

We can characterize weakly distal flows $(X,T)$ as flows which have only singletons as strongly proximal sets and if such flows are not distal then they must admit non trivial proximal sets. Note that for any minimal ideal $I$ in $E(X)$, if for $x \neq y \in X$ we have $ux = uy$ for $u \in J_I$, then for all $p \in I$, $px = pu x = pu y = py$. If this happens for all minimal ideals $I$ then $\{x,y\}$ is a strong proximal set. But for weakly distal flows the strongly proximal sets are just singletons, and so we must have a minimal ideal $I$ for which $ux \neq uy$ for some $u \in J_I$ whenever $x \neq y$. Thus $J$ seperates points in $X$. Further if the flow is weakly distal but not distal then there exists $x \neq y \in X$ with $ux = uy$ for some $u \in J$ and so the elements of $J$ are not all injective hence not invertible.

\bigskip

\section{Examples}

In this section we look into some examples. We observe that it is quite possible that $P \neq SP$ or $SP$ is not closed.

\smallskip

[Note that in all these examples, since both $P$ and $SP$ are reflexive, it is enough to mention $(x,y) \in P (SP)$, since $(y,x) \in P(SP)$ follows. Also, the complement of $  P(SP) $ is $D(WD)$ and so need not be mentioned.]

\bigskip

We first take up a non-minimal example:

\begin{example} \label{CC}
	Consider the circles $ C_n, \ (n = 0, 1, 2, \ldots ) $ defined by
	$ x^2 + (y-i)^2=(1-n/(n^2 + 1))^2  \subset \C$. Let $O_n = i(1-n/(n^2 + 1)) \in \C$ be a point on $C_n$ for all $n$. A point $ Q $ on $ C_n $ will be	given by the coordinates $ ( n, \alpha) $ where $ \alpha $ is the angle  subtended by line $\bar{O_nQ}$ from the horizontal. Thus $C_n = \{(n, \alpha): 0 \leq \alpha < \pi\}$, with $(n,\pi)$ and $(n,0)$ identified.
	
	Consider the circles $ D_n, \ (n =  2, 3, \ldots ) $ defined by
	$ x^2 + (y-i)^2=(1/n)^2  \subset \C$. Notice that $D_2 = C_1$. Let $P_n = i- 1/n \in \C$ be a point on $D_n$ for all $n$. A point $ Q $ on $ D_n $ will be	given by the coordinates $ ( 1/n, \alpha) $ where $ \alpha $ is the angle  subtended by line $\bar{P_nQ}$ from the horizontal. Thus $D_n = \{(1/n, \alpha): 0 \leq \alpha < \pi\}$, with $(1/n,\pi)$ and $(1/n,0)$ identified.
	
	Note that all the circles $C_n$ and $D_n$ are concentric circles with centre at $i$ in $\C$. In fact $\{i\}$ can also be considered as a circle of radius $0$ centred at $i$.
	
	Define $X = (\bigcup \limits_{n=0}^\infty C_n) \cup (\bigcup \limits_{n=2}^\infty D_n) \cup \{i\}$. Then $X$ is a compact subspace of $\C$.
	
	Let $f: X \to X$ be defined as
	$$f(t,\alpha) = \begin{cases}
		(t, \alpha) \ \text{if} \ t=0;\\
		(t+2, \alpha + \sin \alpha/t) \ \text{if} \ t \in 2\Z_+ \setminus \{0\};\\		
		(t-2, \alpha + \sin \alpha/t) \ \text{if} \ t \in 2\Z_+ + 1 \setminus \{1\};\\		
		(2, \alpha + t \sin \alpha)  \ \text{if} \ t =1/3;\\
		(1/(t+2), \alpha + t \sin \alpha)  \ \text{if} \ t \in \{1/n; n = 2,4,6 \ldots \};\\
		(1/(t-2), \alpha + t \sin \alpha)  \ \text{if} \ t \in \{1/n; n = 5,7,9 \ldots \};\\
		i  \ \text{if} \ (t, \alpha) = i.
		
	\end{cases}$$
	
	Note that $f$ is identity on the circle $C_0$ and the point $i$.
	
	Then $f$ is a homeomorphism on $X$ and the only almost periodic points for  the cascade $(X,f)$ are the points on the circle $C_0$ and the point $i$ which are also the fixed points for $f$.
	
	\bigskip	
	
	For the cascade $(X,f)$ we note that:
	
	All points   $ \ \{(t,\alpha): t \in (2\Z_+  \cup \{1/n; n = 3,5,7 \ldots \}) \setminus \{0\}; 0 \leq \alpha < \pi\}$ are forward  asymptotic to $C_0$ and backward asymptotic to $\{i\}$. All points   $ \ \{(t,\alpha): t \in ((2\Z_+ + 1) \cup \{1/n; n = 2,4,6 \ldots \}); 0 \leq \alpha < \pi\}$ are forward  asymptotic to $\{i\}$ and backward asymptotic to $C_0 $. Thus, here the mappings of circles can be described as:
	
	\begin{center}
	
	$ i \longleftrightarrow i $
	
$ 	C_0 \longleftrightarrow C_0 $
	
$ 	C_2 \longrightarrow C_4 \longrightarrow C_6 \longrightarrow C_8 \longrightarrow \ldots $
	
	$ C_1 \longleftarrow C_3 \longleftarrow C_5 \longleftarrow C_7 \longleftarrow C_9 \longleftarrow \ldots $

	$ C_1 = D_2 \longrightarrow D_4 \longrightarrow D_6 \longrightarrow D_8 \longrightarrow \ldots $
	
$ 	C_2 \longleftarrow D_3 \longleftarrow D_5 \longleftarrow D_7 \longleftarrow D_9 \longleftarrow \ldots $
	
	\end{center}

	\smallskip

	\noindent	$ {P = \Delta \ \bigcup \  \{((t_1,\alpha), (t_2,\beta)): t_1, t_2 \in \left(2\Z_+ + 1 \cup \{1/n; n = 2,4,6 \ldots \}\right); 0 \leq \alpha, \beta < \pi\}}$
	
	$\bigcup \  \{((t,\alpha), i): t \in \left(2\Z_+ + 1 \cup \{1/n; n = 2,4,6 \ldots \}\right); 0 \leq \alpha < \pi\}$
	
	$ {\bigcup \ \{((t,\alpha), i):  t \in \left(2\Z_+  \cup \{1/n; n = 3,5,7 \ldots \}\right) \setminus \{0\}; 0 \leq \alpha < \pi\}   }$

	$ {\bigcup \ \{((t_1,\alpha), (t_2,\beta)):  t_1,t_2 \in \left(2\Z_+  \cup \{1/n; n = 3,5,7 \ldots \}\right) \setminus \{0\}; 0 \leq \alpha, \beta < \pi\}   }$.

	\medskip
	
	\noindent And $\  \  SP =  \Delta $.

	\medskip

	We note that  $P$ is neither closed in $X \times X$, nor is $P$  an equivalence relation.	But, $SP$ is trivial.

\smallskip

We note that here $(X,f)$ is weakly distal, but not minimal. Note that all points in  $C_0 \cup \{i\}$ are almost periodic points here.

\bigskip 

Since $P$ is not an equivalence relation here, we note that $E(X)$  will have  many minimal ideals $I_k$, each giving  infinitely many $I_k-$proximal sets since each such $I_k-$proximal set will contain only one  point of the circle $C_0$.

However  the sets $A = \{(t,\alpha): t \in (2\Z_+  \cup \{1/n; n = 3,5,7 \ldots \}) \setminus \{0\}; 0 \leq \alpha < \pi\} \cup \{i\}$ being backward asymptotic to $i$ and $B = \{(t,\alpha): t \in ((2\Z_+ + 1) \cup \{1/n; n = 2,4,6 \ldots \}); 0 \leq \alpha < \pi\} \cup \{i\}$ being forward asymptotic to $i$ are  proximal sets.
 
\end{example}

We now consider an example of a weakly mixing, minimal system.

\begin{example} \label{MT}
	We look into a substitution system,  the square of the Morse-Thue substitution. This is a continuous substitution $Q$ defined by the rule
	$$Q(0) = 0110, \ Q(1) = 1001.$$
	
	We have four bi-infinite sequences that serve as fixed points of $Q$,
	
	\begin{align*}
		a = \ldots 1001.1001 \ldots \\
		b = \ldots 0110.1001 \ldots \\
		\bar{a} = \ldots 0110.0110 \ldots\\
		\bar{b} = \ldots 1001.0110 \ldots \\
	\end{align*}
	
	where $\bar{y}$ denotes the \emph{dual} of $y$ in $\{0,1\}^{\Z}$.
	
	If $x$ denotes any one of the fixed points of $Q$, then $X = \overline{\cO(x)}$ can be defined uniquely for any $x \in \{a,b,\bar{a}, \bar{b}\}$. The system $(X, \sigma)$ is a minimal subsystem of the \emph{$\emph{2}-$shift}.
	
	See the details in \cite{An}, 
	
	Here $E(X)$ has exactly four minimal  idempotents $u_1, v_1, u_2, v_2$ such that they  are the identity off the orbits of $a,b,\bar{a}, \bar{b}$ and on the orbits of $a,b,\bar{a}, \bar{b}$ are defined as:

	\begin{align*}
		u_1: a \to b ; \ \bar{a} \to \bar{b} ; \ b \to b ; \ \bar{b} \to \bar{b}\\
		v_1: a \to \bar{b} ; \ \bar{a} \to b ; \ b \to b ; \ \bar{b} \to \bar{b}\\
		u_2:  a \to a ; \ \bar{a} \to \bar{a} ; \ b \to a ; \ \bar{b} \to \bar{a}\\
		v_2: a \to a ; \ \bar{a} \to \bar{a} ; \ b \to \bar{a} ; \ \bar{b} \to a\\ 
	\end{align*}
	
	Note that,  $E(X)$ here has two minimal left ideals $I, K$ with $u_1, u_2 \in I$ and $v_1, v_2 \in K$. Notice all the points off the orbits of $a,b,\bar{a}, \bar{b}$, which are precisely the distal points in $(X,\sigma)$. Here, $J=\{u_1,u_2,v_1,v_2\}$, and for all $x \in X$ we have $P[x] =Jx$. Hence,

	$$P = \Delta \cup \{(a,b) , (a, \bar{b}),  (b, \bar{a}), (\bar{b}, \bar{a}) \ \text{and their orbits}\}, \text{and} \ SP = \Delta.$$

	Thus, $P$ is neither an equivalence relation nor is closed in $X \times X$. However, for any $u \in J$, there exists a $v \in J$ such that $(vx,vux) \notin \Delta$ and so $SP$ is trivial in $X \times X$.
	
	Note that here for the two minimal ideals $I$ and $K$ we have infinitely many $I-$proximal sets and $K-$proximals sets. 
\end{example}

We note that in the example above, $P$ was discrete and so $WD$  densely contains $D$.

However it is possible that $SP = \Delta$ even when $P$ is a dense subset of $X \times X$, as in the next example.

\begin{example} \label{PL} We take this example from \cite{G}.
	
	The group $ T = SL_2(\R) $  of $ 2 \times 2 $ real matrices with determinant one acts naturally on the space  $ X $ of all rays emanating from the origin in $ \R^2 $ (considered as the space of all $ 1 \times 2$ real row vectors with the action being   the matrix multiplication with the row vector). This action of $ T $	on the compact space $ X $  is  minimal.
	
	$ X $ can be identified with the unit circle in $ \R^2 $, and it can be shown that the minimal ideals	of $ E(X) $ are in one-to-one correspondence with the partitions of $ X $ into two
	complementary half closed arcs in the following way - If $ I $ and $ I' $ are two such	complementary half closed arcs then the elements of the corresponding	minimal ideal $ M  $ are in one-to-one correspondence with the points of $ X $. If $ x \in I $
	then the corresponding element $ v_x \in M $ maps all of $ I $ onto $ x $ and all of $ I' $ onto its
	antipodal point $\breve{x}$(thus  $ v_x $ is an idempotent). 
	
	Here $J_M = \{v_{x_{[I;I']}}: x \in X, \ X = I \cup I'\}$ with each minimal ideal $M$ corresponding to the partition $X =I \cup I'$, and $J = \bigcup \limits_{X = I \cup I'} J_M$.	Now 
	$$P = \bigcup \limits_{v \in J} \ \{(x,vx):  \ x \in X\} = X \times X \setminus \{(x, \breve{x}): x \in X\}.$$ 
	
	And for $x \neq y \in X$, there exists $v \in J_{\hat{M}}$ with $\hat{M}$ corresponding to some partition of $X$ such that $v(x,y) = (x,y)$. This gives,
	
	$$SP = \Delta.$$
	
	Note that $P$ is dense here, but $SP = \Delta$.
	
	The $I-$proximal sets here can be given by all  partitions of $X$ into two complementary half closed arcs that correspond to the minimal ideals. 
	
\end{example}

We now give an  example of  minimal flow with $P = SP$.

\begin{example} \label{F}
	We take the example of a ``non-homogeneous" minimal set, which is a
	modification of a construction of E. E. Floyd as given in \cite{AUS}.
	
	If $B=[0,2] \times [0,1] \subset \R^2$ , and let $ h ( B ) $
	denote a union of $ 3 $ disjoint rectangles $ h ( B ) = B_0 \cup B_1 \cup B_2 $, where $	B_0 = [0,2/5] \times [0,1/2], \ 	B_2 = [4/5,6/5] \times [0,1], \ B_3 = [8/5, 2] \times [1/2,1]$. That is, $  h ( B ) $ is obtained from $ B $ by 	deleting the ``middle fifths" rectangles and then deleting the top and	bottom halves respectively of the remaining left and right rectangles.
	
	Now let $ B^{(O)} = B  $, $B^{(1)} = h(B^{(0)})$ and define inductively $ 	B^{(n+l)} = h(B^{(n)})$. 
	
	Thus $ B^{(n)} $ consists of $ 3^n $ disjoint rectangles
	$B_k^{(n)}$ $(k = 0,1,2,.. .,3^n -1)$ arranged from left to right as:
	
	$B^{(1)} = B_0^{(1)} \cup B_1^{(1)} \cup B_2^{(1)} = h(B^{(0)}) = B_0 \cup B_1 \cup B_2$.
	
	$h(B_0^{(1)}) = B_0^{(2)} \cup B_3^{(2)} \cup B_6^{(2)} $, $h(B_1^{(1)}) = B_1^{(2)} \cup B_4^{(2)} \cup B_7^{(2)} $, and $h(B_2^{(1)}) = B_2^{(2)} \cup B_5^{(2)} \cup B_8^{(2)} $.
	
	$ \ldots $
	
	In general, $h(B_j^{(n)}) = B_j^{(n+l)} \cup  B_{j+3^n}^{(n+l)} \cup B_{j+2 3^n}^{(n+l)}$ from left to right.
	
	$ X = \bigcap \limits_{n=0,1,2, \ldots} B^{(n)} $ is a compact metric space, which consists of vertical
	line segments, and	degenerate (single) points. The metric space $ X $ is nonhomogeneous 	-- there is one segment of length 1 and, for
	every $ n > 0 $, infinitely many of length $ 1/2^n $, and the rest are singletons.	
	
	Define a homeomorphism $ T $ of $ X $ as  by permuting the rectangles
	$ B_j^{(n)} $. These permututations of the rectangles
	$ 	B_j^{(n)} \to  B_{j+1}^{(n)}  $ (where $ j+l $ is considered $ \mod 3^n $ ) induce a map of X to itself. The cascade (X,T) is minimal. Here,

	$$P = \{ (x , x’)  \in X \times X:  x , x' \in X \ \text{are on the same vertical segment}\}.$$
	
	$P$ is an equivalence relation and so $SP = P \neq \Delta$ here. 
	
Note that the vertical lines give the proximal sets which are also the strongly proximal sets.	
\end{example}

Next we give an example of a minimal flow with $SP$ non-trivial, yet closed.

\begin{example} \label{Shift}
	We take this example from \cite{EG}.
	
	Let $ X = \{0, 1, 2\}^\Z $ and $ \sigma $  denote the shift on $ X $. The symmetric	group $ \cS_3 $ acts on $ X $ by permuting values of the zero coordinate of  points $ x \in  X $. Let  $ \tau_{n,m} $ be the homeomorphism of permutating the $ nth $ and $ mth $ coordinate of elements of $ 	X $. 
	
	Let $  G = \langle \sigma , \cS_3, \tau_{n,m} : n, m \in \Z \rangle $ be the subgroup of homeomorphisms of $ X $
	generated by $\sigma, \cS_3 $ and all $\tau_{n,m}$.

	Every homeomorphism $ g \in G $ is characterized by the property,
	$$ \forall \ x, y \in X [x_n \neq y_n \ \forall \ n \in \Z \Rightarrow (gx)_n \neq (gy)_n \ \forall \ n \in \Z].$$

	We follow these terminologies: call
	a pair of points $ (x, y) \in X \times X  $ \emph{an edge} if $ x_n \neq y_n $ for every $ n \in \Z $; if $ x_n \neq y_n $ for
	infinitely many $ n \in Z  $, then $ x $ and $ y $ are called \emph{opposed}; and if $ x $ and $ y $ are not opposed, i.e. if $ x_n = y_n $ for all but 	finitely many $ n \in Z  $, they are called \emph{agreeable}.
	
	\bigskip
	
	As is proved in \cite{EG}, $ ( X, G )$ is minimal.
	
	\bigskip
	
	Let $ \bar{0}, \bar{1} $ and $ \bar{2} $ denote the points of $ X $ all whose coordinates are $  0, 1 $	and $ 2 $ respectively.

	\bigskip
	
	As proved in \cite{EG}: Given  points $ \{\bar{0}=x_0, x_1, x_2	, x_3\} $ in $ X $, there exists a sequence $ 	g_n \in G  $ with $ \lim g_nx_j \in \{\bar{0}, \bar{1}, \bar{2}\} $ for $ j = 0, 1, 2, 3 $.
	
	We recall the proof mentioned there. If all of $x_0,x_1,x_2,x_3$ are agreeable then $ \lim \limits_{k \to \infty} \sigma^kx_j = \bar{0} $ for $ j = 0, 1, 2, 3 $. 	If there exists $ j $ with $ x_j $
	opposed to $ \bar{0} $ (i.e. $ x_j $ has infinitely many coordinates different from $ 0 $), we can
	apply a sequence of $ \tau $ -permutations, elements of $ \cS_3 $ and various powers of $\sigma$ to
	the pair $ \{\bar{0}, x_j\} $ to have the limit in $ \{\bar{0}, \bar{1}\} $.  In either
	case passing to a further subsequence we now have:
	$ { \lim \limits_{n \to \infty} g_nx_j	} \subset \{\bar{0}, \bar{1}, y_2, y_3\} $	for some $ y_2	, y_3 \in  X $. If either $ y_2 $ or $ y_3 $ has infinitely many coordinates with the
	value $ 2 $, we can similarly pass to a limit which is a subset of $  \{\bar{0}, \bar{1}, \bar{2}, z_3\} $ for some $ z_3 \in X $. Otherwise we can pass to a subset of $ \{\bar{0}, \bar{1}\} $. Finally, from $ \{\bar{0}, \bar{1}, \bar{2}, z_3\} $ 	we can get into $ \{\bar{0}, \bar{1}, \bar{2}\} $. 
	
	\bigskip
	
	Further discussions in \cite{EG} based on the above discussion  proves that $\Omega(X)	 = $ all edges in $X \times X$.
	
	\bigskip
	
	Note that there exists  $z= \ldots 00011000.110001110001111000 \ldots \in X$ that is agreeable to $\bar{0}$ and opposed to $\bar{1}$. Then by the discussion above  $\lim \limits_{n \to \infty} \sigma^n z = \bar{0}$. And a sequence $\{g_n\}$ of $ \tau $ -permutations and various powers of $\sigma$ to have $\lim \limits_{n \to \infty} g_n z = \bar{1}$. Hence $(\bar{0},z) \in P(X)$ and $(z, \bar{1}) \in P(X)$. But $(\bar{0},\bar{1})  \notin P(X)$. Hence, $P$ is not an equivalence relation here. The same argument also implies that $P$ is not closed.
	
	\bigskip
	
	Here, $SP(X) \neq P(X)$.
	
	Let $(x,y) \in X \times X \setminus \Omega(X)$ be such that both $x,y$ are agreeable to $ \bar{0} $. Then $\overline{G(x,y)} \cap \Omega(X) \ni (\bar{0},\bar{0}), (\bar{1},\bar{1})$ or $(\bar{2},\bar{2})$. Thus, $(x,y) \in SP(X)$.
	
	If $x$ is agreeable to $\bar{0}$ and $y$ is opposed to $\bar{0}$, but agreeable to either $\bar{1}$ or $\bar{2}$. then again by the discussion above we have $\overline{G(x,y)} \cap \Omega(X) \ni (\bar{0},\bar{1})$ or $(\bar{0},\bar{2})$. Thus, $(x,y) \notin SP(X)$.
	
	If $x$ is agreeable to $\bar{0}$ and $y$ is opposed to all $\bar{0}, \bar{1}, \bar{2}$, then by the discussion above there exists a sequence $\{g''_n\}$ of $ \tau $ -permutations, elements of $ \cS_3 $ and various powers of $\sigma$ to have $\lim \limits_{n \to \infty} g''_n y = \bar{1}$( or $\bar{2}$). Then, $(x,y)
	\notin SP(X)$.

	If $x,y$ are both opposed to all $\bar{0}, \bar{1}, \bar{2}$, but are mutually agreeable then $\overline{G(x,y)} \cap \Omega(X) \ni (\bar{0},\bar{0}), (\bar{1},\bar{1})$ or $(\bar{2},\bar{2})$. Thus, $(x,y) \in SP(X)$.

	If $x,y$ are both opposed to all $\bar{0}, \bar{1}, \bar{2}$, and are also mutually opposed  then there are infinitely many mutually distinct coordinates where $x$ and $y$ both take the value say $1$. Then by the discussion above there exists a sequence $\{g_n\}$ in  $ G $  to have $\lim \limits_{n \to \infty} g_n y = \bar{1}$, but no subsequence of $\{g_n\}$ will have $\lim \limits_{k \to \infty} g_{n_k} x = \bar{1}$. Thus $(x,y) \notin SP(X)$.
	
	Thus, $SP(X) = \{(x,y) \in X \times X: x,y \ \text{are mutually agreeable}.\}$
	
	We now prove that $SP(X)$ is closed in $X \times X$. Let $(x,y) \in \overline{SP(X)}$ then there exists a sequence $\{(x_n,y_n)\}$ in $SP(X)$ such that $(x_n,y_n) \to (x,y)$. Now for each $n \in \N$, $x_n$ and $y_n$ are mutually agreeable. Since $x_n \to x$ and $y_n \to y$, $x$ and $y$ are also mutually agreeable. Thus $(x,y) \in SP(X)$ and so $SP(X)$ is closed.
	
	Thus, here $P$ is neither an equivalence relation nor is closed in $X \times X$. However, $SP$ is closed in $X \times X$.
	
	Note that here the strongly proximal sets comprise of those sequences that are mutually agreeable. The proximal sets are more complicated.
	
\end{example}

Finally, we consider a classical weakly mixing, minimal transformation.

\begin{example} \label{Chacon}
	We consider the symbolic version of the classical \emph{Chac\'{o}n's transformation}. We refer to \cite{dJK, VER} for details of this example.
	
	\noindent This example  is a subshift of the $2-$shift $ (\{0, 1\}^{\Z}, \sigma) $
	with  blocks of $0's$ spaced systematically by a $1$.
	
	 We define 	finite blocks $ B_0, B_1, B_2, \ldots $ inductively by setting
	$ B_0 = 0, B_1 = 0 0 1 0, B_{k+1} = B_k B_k 1 B_k $ for all $k \in \N$. Each $B_k$ is called a $k-$block. And we note that $B_{k+1}$ comprises of three $k-$blocks with $1$ spacing the second and third occurences of $B_k$.
	
\centerline{Define	$X:=\{x \in \{0,1\}^{\Z}:  x_{[i,j]} $, for  all $i < j \in \Z$, is a subword of $ B_k $ for	some $k \in \N \}.$}
	
\noindent Then	$ X $ is clearly a closed invariant subset of the $2-$shift.
	
	 We explicitly describe $X$. Since each $B_{k+1}$ begins with $B_k$ there exists a unique right sequence $B_{\infty}$ such that $B_k$ is the  $|B_k|-$ prefix of $B_{\infty}$. Similarly, each $B_{k+1}$ ends with $B_k$ and so there exists a unique left sequence $B_{-\infty}$ such that $B_k$ is the  $|B_k|-$ suffix of $B_{-\infty}$.
	 
	 Again, we refer to \cite{dJK, VER} for the proof of the observation that $B_k$'s do not overlap and in a spaced concatenation of $ k- $blocks, each block occur only at its natural	 positions with the spacer occuring systematically. 
	
	We take $\mathbf{x_1} = B_{-\infty}B_\infty$ and $\mathbf{x_2} = B_{-\infty}1B_\infty$. Clearly $\mathbf{x_1, x_2} \in X$. 
	
	We can uniquely define all $x \in X$. 	Let $\xi \in \{1,2,3\}^{\N}$.
	
	 Define $D_0 := B_0 = 0$ and inductively $D_{k+1} := \begin{cases} D_kB_k1B_k \ \text{if} \ \xi_k = 1\\
	B_kD_k1B_k \ \text{if} \ \xi_k = 2\\ B_kB_k1D_k \ \text{if} \ \xi_k = 3 \end{cases}$ for all $k \in \N$.

Note that $D_k = B_k$ for all $k \in \N$. We construct $\xi^*$ corresponding to each $\xi$ starting from $D_0$.

1. If $\xi_k = 1$ for all $k \geq n_0$ for some $n_0 \in \N$ then $\xi^* = B_\infty$ as the sequence grows only to the left. This gives $x \in X$ with $x_{[n_0,\infty)} = B_\infty$. And so by the systematic concatenation of the blocks, we have $x = \sigma^n(\mathbf{x_1})$ for some $n \in \Z$.

2.  If $\xi_k = 3$ for all $k \geq m_0$ for some $m_0 \in \N$ then $\xi^* = B_{-\infty}$ as the sequence grows only to the right. This gives $x \in X$ with $x_{(-\infty, m_0]} = B_{-\infty}$. And so by the systematic concatenation of the blocks, we have $x = \sigma^m(\mathbf{x_1})$ or $x = \sigma^m(\mathbf{x_2})$ for some $m \in \Z$.

3. For all other $\xi$, we note that for a defined $D_k$, it can be extended to $D_{k+1}$ uniquely depending on $\xi_{k+1}$. Thus the resulting $\xi^*$ is extended both to the left and right. This gives $x \in X$ with $x = \sigma^l(\xi^*)$ for some $l \in \Z$.

Thus $X$ consists precisely all the orbits of $\xi^*$ when $\xi$ is not eventually $1$ or $3$ together with the orbits of $\mathbf{x_1}$ and $\mathbf{x_2}$. Trivially this $X$ is minimal.

\bigskip

Now for $x,y \in X$, with $x = \sigma^p(y)$ for some $p \in \N$, we note that $(x,y) \in D$ and so $(x,y) \notin SP$. We again refer to \cite{dJK, VER} for the detailed proof that for $x,y \in X$ with one of them not equal to  $\mathbf{x_1}$ or $\mathbf{x_2}$ and both not sharing the same orbit, we have $\overline{\bigcup \limits_{n=1}^\infty \ (\sigma \times \sigma)^n(x,y)} = X \times X$. And so  $(x,y) \notin SP$. And finally for $x,y \in X$ with $\sigma^n(x) = \mathbf{x_1}$ and $\sigma^m(y) = \mathbf{x_2}$ for some $n,m \in \Z$, we observe that for arbitrary large $k > 0$, $\sigma^{k+n}(x)$ and $\sigma^{k+m}(y)$ share the same central block. Thus, there exists $\{k_i\}$ with $k_i \nearrow \infty$ such that $\sigma^{k_i + n}(x) \rightarrow z \leftarrow \sigma^{k_i +m}(y)$. If $n \neq m$ then clearly $(x,y) \notin  SP$. If $n = m$, we note that for arbitrary large $k > 0$, $\sigma^{k+n}(x)$ and $\sigma^{k+n+1}(y)$ share the same central block, and so as discussed above $(x,y) \notin  SP$.

Thus $SP = \Delta$, and so $(X,\sigma)$ is weakly distal.

Note that here since $(X,\sigma)$ is weakly mixing and minimal so $P(X)$ is dense in $X \times X$.

\end{example}
\vspace{12pt}
\bibliography{xbib}

\end{document}